\documentclass[preprint,12pt,number]{elsarticle}
\usepackage{setspace}
\usepackage[
    a4paper,
    left=2cm,
    right=2cm,
    top=2.5cm,
    bottom=2.5cm
]{geometry}

\usepackage{graphicx}%
\usepackage{multirow}%
\usepackage{subcaption}
\usepackage{amsmath,amssymb,amsfonts}%
\usepackage{amsthm}%
\usepackage{mathrsfs}%
\usepackage[title]{appendix}%
\usepackage{xcolor}%
\usepackage{textcomp}%
\usepackage{manyfoot}%
\usepackage{booktabs}%
\usepackage{algorithm}%
\usepackage{algpseudocode}%
\usepackage{listings}%

\usepackage{pgfplots}
\usepgfplotslibrary{groupplots,statistics}
\usepackage{pgfplotstable}
\pgfplotsset{compat=1.18}
\usepgfplotslibrary{fillbetween}
\definecolor{colOut}{RGB}{31,119,180}
\definecolor{colEmg}{RGB}{255,127,14}
\definecolor{colInp}{RGB}{44,160,44}

\def\medoutexam{15.00}

\def\medoutrep{55.00}

\def\medemgexam{8.00}

\def\medemgrep{20.00}

\def\medinpexam{10.00}

\def\medinprep{25.00}

\usepackage{hyperref}
\usepackage{xurl}
\usepackage{tikz}
\usetikzlibrary{arrows.meta, shapes.geometric, positioning, patterns, calc}
\usepackage{booktabs}

\usepackage{xcolor}
\definecolor{softgreen}{RGB}{0,100,0}

\begin{document}

\begin{frontmatter}



\title{A Predictive-Prescriptive Analytics Framework for Fair Computed Tomography Scheduling and Radiologist Workload Allocation} 


\author[p]{Ludovico Ambrosi}
\author[p]{Chandra Bortolotto}
\author[u]{Sara Cambiaghi\corref{cor1}}
\cortext[cor1]{Corresponding author}
\ead{sara.cambiaghi01@universitadipavia.it}
\author[p]{Luisa Carone}
\author[u]{Davide Duma}
\author[p]{Lorenzo Preda}

\affiliation[p]{organization={Department of Radiology, Fondazione I.R.C.C.S. Policlinico 
San Matteo},
            addressline={Viale Camillo Golgi, 19}, 
            city={Pavia},
            postcode={27100}, 
            country={Italy}}
\affiliation[u]{organization={Department of Mathematics “F. Casorati”, University of Pavia},
            addressline={Via Ferrata, 5}, 
            city={Pavia},
            postcode={27100}, 
            country={Italy}}

\begin{abstract}
Scheduling follow-up Computed Tomography (CT) examinations requires balancing
two competing objectives: assigning patients as close as possible to their
recommended examination dates while ensuring an equitable distribution of
radiologists' workload. Existing approaches optimize scanner utilization or
patient waiting times, overlooking reporting activities and the need to balance
fairness across multiple stakeholders. This paper proposes a
predictive-prescriptive framework for fairness-aware follow-up CT scheduling.
Patient-specific Machine Learning (ML) models are first developed to predict
both examination and reporting durations. These predictions are then embedded
into a multi-objective Mixed-Integer Linear Programming (MILP) model that
simultaneously minimizes deviations from patients' preferred examination dates
and balances radiologists' reporting workloads through a lexicographic min-max
fairness criterion. We derive a dominance reduction property within an
$\varepsilon$-constraint framework that substantially reduces the number of
optimization problems required to generate the Pareto frontier. Computational
experiments based on data from a real-world emergency radiology department show
that allowing patients a scheduling flexibility of only one to two days is
sufficient to substantially improve workload equity among radiologists while
preserving timely access to follow-up examinations. The proposed dominance
reduction strategy eliminates most $\varepsilon$-constraint evaluations without
affecting the Pareto frontier. Finally, evaluating predictive models through
downstream optimization regret demonstrates that XGBoost provides the most
effective support for scheduling decisions, outperforming models that achieve
lower prediction errors according to conventional predictive metrics.
\end{abstract}

\begin{highlights}
\item A fairness-aware framework jointly optimizes scheduling and workload allocation.
\item Machine learning predicts tasks durations, integrated into the optimization model.
\item A dominance reduction property improves computational efficiency.
\item The proposed framework is validated using real-world emergency radiology data.
\end{highlights}

\begin{keyword}
Predictive analytics \sep Machine learning \sep Prescriptive analytics \sep Healthcare scheduling \sep Computed tomography \sep Mixed-Integer linear programming
\end{keyword}

\end{frontmatter}

\section{Introduction}
\label{sec:FAIR:intro}

Computed Tomography (CT) is a fundamental component of emergency radiology, providing rapid diagnostic support for trauma, stroke, oncological emergencies, and other time-critical conditions. In tertiary hospitals, CT scanner capacity is primarily devoted to Emergency Department (ED) patients, whose unpredictable arrivals require immediate access to imaging resources \cite{DanLantsman2022}. Alongside this unscheduled demand, radiology departments also have to accommodate examinations requested for hospitalized patients \cite{MacDonald2013} and outpatient follow-up CT examinations \cite{Kavandi2023}.

A considerable proportion of outpatient CT examinations performed within emergency radiology departments are scheduled as follow-up investigations after a previous ED visit or hospital admission. Similarly, some inpatients may require scheduled follow-up examinations during their hospital stay. These examinations represent a continuation of the diagnostic pathway initiated during the acute episode rather than routine elective imaging. Follow-up CT may be required to assess the evolution of findings, confirm or exclude suspected diagnoses, evaluate treatment response, or monitor lesions whose interpretation depends on temporal changes. Consequently, radiologists or referring physicians typically recommend that the examination be performed after a specified interval \cite{Kavandi2023}. This recommendation naturally defines a \textit{target examination date}, representing the clinically preferred timing for the follow-up and, whenever feasible, accommodating inpatients and outpatient scheduling preferences \cite{Doshi2023}. 

Scheduling follow-up CT examinations within an emergency radiology department is operationally challenging because available capacity is neither fixed nor predictable. The workload of CT scanners is strongly influenced by the fluctuating demand generated by emergency patients, whose arrivals vary across days of the week and work shifts \cite{Cambiaghi2026}. Consequently, the capacity that can be allocated to follow-up examinations changes dynamically over time. This complexity is further amplified by the substantial variability in both examination durations and, even more importantly, radiologist reporting times, which depend on patient-specific clinical conditions and imaging findings. Balancing these heterogeneous and uncertain workloads while assigning inpatients and outpatient follow-up appointments close to their target dates therefore represents a challenging operational problem, requiring a trade-off between clinical timeliness and efficient utilization of radiology resources \cite{Aloini}.
For simplicity, throughout the remainder of the paper, we refer to both inpatients and outpatients requiring follow-up examinations as outpatients. Conversely, the term inpatients refers to hospitalized patients requiring urgent care whose examinations cannot be scheduled in advance.

These operational challenges have motivated extensive research on radiology scheduling and healthcare resource allocation. Appointment scheduling and resource allocation in radiology have been extensively investigated using operations research techniques, including queueing theory, stochastic optimization, simulation, and Markov decision processes \cite{Patrick, Salemi, Zhang}. These approaches have successfully improved resource utilization and reduced waiting times under uncertainty. In particular, integrated optimization models have been proposed to jointly coordinate patient appointments and healthcare resources, such as technologist schedules, demonstrating the benefits of considering multiple stakeholders simultaneously \cite{Bentayeb2023}.

Nevertheless, most existing scheduling models optimize scanner utilization or waiting times without explicitly considering patient-specific target follow-up dates, and typically model scanner capacity as the primary limiting resource, overlooking radiologist reporting workload as a potential operational bottleneck \cite{Patrick, Salemi, Zhang}. 

Recent advances in Machine Learning (ML) have enabled accurate prediction of operational parameters in healthcare, including surgical durations \cite{Ala}, hospital length of stay \cite{Huang}, appointment cancellations \cite{Rothenberg}, and CT examination durations \cite{Bhattacharjee, Wang2024}. However, comparatively little attention has been devoted to predicting radiologist reporting times \cite{Cowan2013}, and these predictive models have rarely been integrated into optimization frameworks for radiology appointment scheduling. 
Consequently, an integrated predictive-prescriptive framework that jointly accounts for patient-specific target follow-up dates, predicted examination and reporting workloads, and stakeholder-specific fairness objectives is still lacking.

Reliable estimates of examination and reporting times are essential for effective operational decision-making, but prediction alone is insufficient. These estimates must be embedded within scheduling models capable of balancing clinical and organizational requirements. In the context of follow-up CT scheduling, this gives rise to two conflicting objectives: assigning examinations as close as possible to their target dates while promoting an equitable distribution of reporting workloads across radiologists and working shifts. 

Motivated by these considerations, this study proposes an integrated predictive-prescriptive framework for follow-up CT scheduling in emergency radiology. ML models are first used to predict patient-specific CT examination and reporting durations. The resulting estimates are then incorporated into a multi-objective Mixed-Integer Linear Programming (MILP) model that explicitly addresses multi-stakeholder fairness by jointly optimizing patient- and radiologist-oriented objectives. Specifically, the model minimizes deviations from patients' target examination dates while promoting an equitable distribution of reporting workloads across radiologists and working shifts through a lexicographic min–max criterion \cite{Bertsimas2011}.

Simultaneously optimizing these two stakeholder-specific notions of fairness gives rise to a challenging multi-objective optimization problem, requiring the efficient computation of non-dominated solutions that characterize the trade-offs between timely patient access to follow-up CT examinations and equitable distribution of reporting workloads. To efficiently generate the Pareto frontier, we introduce a dominance reduction property that significantly reduces the search space explored by the $\epsilon$-constraint method. The proposed framework provides a practical predictive-prescriptive decision-support tool that integrates predictive analytics with mathematical optimization for fair and efficient follow-up CT scheduling in emergency radiology.

The remainder of the paper is structured as follows. Section \ref{sec:FAIR:literature} reviews the literature most relevant to the present study and identifies the research gaps addressed by this work. Section \ref{sec:FAIR:prediction} presents the predictive models adopted in the study, together with the evaluation metrics used to assess their performance. Section \ref{sec:FAIR:opt} formulates the outpatient scheduling problem and presents the fairness-aware optimization model, the $\epsilon$-constraint reformulation, the dominance reduction property, and the baseline fairness-unaware model. Section \ref{sec:FAIR:results} introduces a case study, presents the results of the integrated predictive-prescriptive framework, and discusses the resulting managerial insights. Finally, Section \ref{sec:FAIR:conclusion} summarizes the main conclusions and outlines directions for future research.

\section{Literature Review}
\label{sec:FAIR:literature}

This section reviews the literature most relevant to this study. We first discuss the application of Artificial Intelligence (AI) and ML to CT workflow management, with a particular focus on patient scheduling. We then review fairness-aware optimization approaches in radiology. Finally, we identify the existing research gap and highlight the contributions of the proposed framework.

\subsection{Artificial Intelligence for Computed Tomography}

The use of AI in healthcare has become increasingly common, including in operating room scheduling \cite{Bellini, Bellini2}, as well as in predicting patient waiting times and appointment delays \cite{Curtis} and prolonged length of stay \cite{Zeleke}.
According to Pierre et al. \cite{Pierre2023}, AI has numerous impactful applications in radiology throughout the imaging workflow. These include applications before image acquisition, such as scan ordering, patient screening, protocol selection, and patient scheduling; during image acquisition, such as patient positioning, intravenous contrast dosing, scan optimization, dose reduction, image post-processing, and workflow support; and during image interpretation, including visualization, quantification, AI-assisted segmentation, radiology reporting, and urgent finding notification \cite{Ferdous2024}.


In this study, we focus on ML applications for patient scheduling. Specifically, ML is used to predict uncertain parameters that are subsequently incorporated into a mathematical optimization model, enabling more accurate scheduling decisions under uncertainty. To this purpose, ML has been applied to predict service times \cite{Daldossi2026, Golmohammadi2023, Lai2024}, identify patients at high risk of no-shows \cite{Srinivas2018, Han2024}, forecast patient demand \cite{Luo2017, Lin2019}, and estimate the performance of appointment schedules \cite{Wang2021, Huang2015}. However, its application to CT scheduling remains limited.

Only a few studies have investigated the use of ML for radiology appointment scheduling. Bhattacharjee et al. \cite{Bhattacharjee} integrated Decision Tree (DT) predictions into a discrete-event simulation to evaluate scheduling policies for walk-in and emergency patients. Huang et al. \cite{Huang} used DTs to optimize radiology appointment slot lengths based on patient characteristics. 
Wang et al. \cite{Wang2024} used a Random Forest (RF) regressor to predict the duration of individual abdominal Magnetic Resonance Imaging (MRI) examinations, highlighting discrepancies between planned and actual slot lengths. 

This work builds upon the predictive framework introduced in Cambiaghi et al. \cite{Cambiaghi2026}, who compared several machine learning techniques integrated into a mathematical model for outpatient scheduling, while explicitly accounting for emergency and inpatient arrivals. To the best of our knowledge, it is the only study that systematically compares multiple ML techniques for predicting outpatient CT examination durations and integrates the resulting predictions into an optimization model for CT scheduling.
Indeed, most existing studies either focus exclusively on the predictive component or incorporate predictions into simulation-based frameworks. Furthermore, the reporting process is almost always overlooked, and fairness considerations are generally neglected. These gaps highlight the need for approaches that jointly optimize both the examination and reporting phases while incorporating fairness criteria for both patients and radiologists, thereby enabling more comprehensive and data-driven scheduling decisions in radiology departments.

\subsection{Fairness-Aware Optimization in Healthcare Management Application}

Fairness has emerged as an important objective in healthcare scheduling and resource allocation, where decision-makers must balance operational efficiency with the equitable allocation of limited clinical resources \cite{Aringhieri2022, Grot2022, Chorfi2026}. In this context, fairness refers to the equitable allocation of healthcare resources according to transparent and consistent decision rules. It encompasses both distributive justice, which concerns the fairness of allocation outcomes, and procedural justice, which concerns the fairness of the processes through which allocation decisions are made \cite{Guindo2012}.

In healthcare appointment scheduling, fairness has been incorporated into optimization models through different formulations. For example, Yan et al. \cite{Yan2015} proposed a sequential appointment scheduling model that jointly considers patient preferences and service fairness. Their model maximizes the clinic's expected profit while enforcing fairness through a constraint that limits the maximum difference in average waiting times across appointment intervals, thereby promoting a more uniform quality of service while explicitly analyzing the trade-off between operational efficiency and service equity.

Nabavizadeh et al. \cite{Nabavizadeh2024} proposed a stochastic mixed integer programming model for home healthcare routing and scheduling under uncertainty. Besides minimizing operational costs, the model incorporates workload balancing, work regulations, and caregiver satisfaction to promote a more equitable distribution of workload among healthcare personnel, while explicitly accounting for uncertain patient demand and service times through a two-stage stochastic programming framework.

Within operating room scheduling, Kayvanfar et al. \cite{Kayvanfar2025} proposed an integrated framework that combines capacity planning with fair and resilient surgery scheduling. In the capacity-planning stage, a Markovian queueing model determines the required number of operating rooms. Subsequently, a goal-programming model promotes the equitable assignment of patients among surgeons, while reducing surgeons’ overtime and idle time and accounting for cancellation risks through the identification of backup surgeons.

In the context of diagnostic imaging, multiple resources including CT scanners, radiologists, technicians, and other healthcare personnel must be allocated to patients, who may also have different preferences and priorities. Consequently, fairness should be considered alongside efficiency when designing scheduling and resource allocation policies. 
Zhou et al. \cite{Zhou2015} proposed an integer linear programming model to allocate MRI capacity among different patient groups, classified according to the body part being examined. The proposed model minimizes the total cost, including patient waiting costs, patient rejection costs, and machine idle costs, while ensuring fairness in the provision of medical services. The authors assume that all patients within the same group have identical, deterministic examination durations. Fairness is enforced by constraining both the deviation of each patient group's rejection rate from the average rejection rate and the deviation of each machine's workload from the average workload to remain below predefined thresholds.

Overall, the existing literature addresses fairness in healthcare applications primarily through single-stakeholder formulations, emphasizing either equitable patient service or equitable workload allocation among healthcare professionals. Multi-stakeholder fairness has received little attention, particularly in the context of follow-up CT scheduling.

\subsection{Research Gap and Contributions}

To the best of our knowledge, no previous study has integrated ML-based prediction, fairness-aware optimization, and outpatient CT scheduling within a unified decision-support framework.

This work addresses several important gaps in the literature. First, only a limited number of studies have investigated the use of ML techniques to support CT scheduling decisions. Second, existing studies primarily focus on predicting examination durations, whereas reporting times have received little attention despite their significant impact on radiology workflows. Third, most existing approaches rely on relatively simple predictive models and integrate predictions into simulation-based rather than optimization-based decision-support tools. Finally, fairness considerations are rarely incorporated into appointment scheduling models, particularly with respect to both patient allocation and radiologist workload distribution.

The main contributions of this work are summarized as follows:
\begin{itemize}
\item compare multiple ML techniques for predicting both CT examination durations and reporting times, providing a comprehensive assessment of their predictive accuracy and downstream operational performance through a regret-based evaluation;
\item formulate a fairness-aware MILP model that explicitly incorporates both the examination and reporting phases of the CT workflow. The proposed model jointly addresses patient scheduling and radiologist workload allocation while taking into account medically preferred examination days and workload fairness;
\item derive a dominance reduction property that accelerates Pareto frontier generation within the $\epsilon$-constraint method, improving the computational efficiency of the proposed multi-objective optimization model;
\item validate the proposed framework using real-world data from an emergency radiology unit of a large trauma center in Italy and analyze its performance under different demand and case-mix scenarios.
\end{itemize}

\section{Predictive Modeling Framework}
\label{sec:FAIR:prediction}

This section presents the predictive models adopted in the study, together with the evaluation metrics used to assess their performance.

\subsection{Predictive models}
The objective of the predictive component is to estimate, for each patient–radiologist pair, both the CT examination duration and the reporting duration. These predictions constitute the input parameters of the optimization model described in the following section.

Prediction of execution and reporting durations was framed as a regression problem. We trained and compared four ML algorithms and one deep learning model. Specifically, we considered Decision Trees (DT), Random Forests (RF), eXtreme Gradient Boosting (XGBoost), Categorical Boosting (CatBoost), and a Multilayer Perceptron (MLP).
The performance of the model was compared against a baseline that predicts the historical average duration of previous examinations.
In addition, feature importance analysis was performed to identify the predictors that contributed most strongly to duration estimates.


\subsection{Predictive Performance through Regret and Counterfactual Evaluation}
To evaluate the practical impact of prediction errors on scheduling decisions, we assess each ML model using a regret metric within the optimization framework.
Let $obj_{\text{ML}}$ denote the objective function value obtained when the MILP model is solved using the predictions generated by a given ML technique and hyperparameter configuration, with the resulting schedule evaluated against the true task durations. Let $obj_{\text{best}}$ denote the deterministic oracle objective value obtained by solving the MILP directly with the true task durations. The regret is then defined as follows

\begin{equation}
\label{eq:FAIR:regret}
\text{Regret} = \frac{obj_{\text{ML}} - obj_{\text{best}}}{obj_{\text{best}}}.
\end{equation}

The computation of this metric requires a counterfactual evaluation procedure. In the historical dataset, we only observe the actual exam and reporting durations for the specific radiologist assigned to a patient; the durations for alternative, unassigned radiologists remain unobserved. To overcome this limitation and avoid evaluation bias, we we construct a synthetic ground-truth dataset that provides examination and reporting durations for every patient-radiologist pair independently of the ML models' predictions.

For a given patient-radiologist pair, if the assignment matches the historical record, the synthetic duration is set to the true observed duration. Otherwise, we identify a pool of historical patients who were assigned to that specific radiologist and share identical characteristics across the following features:
\begin{itemize}
    \item number of exams;
    \item main exam type (e.g., abdominal CT, intracranial CT angiography, brain CT, chest CT);
    \item department: the hospital department from which the patient was referred;
    \item patient origin: an attribute detailing if the patient arrives from the emergency department (ED), accesses
    the CT scan during hospital stay (inpatient), or for a post-hospitalization check (outpatient).
\end{itemize}
If no exact match exists in the dataset, we hierarchically relax the matching constraints by removing the least critical features sequentially (starting from patient origin) until at least one matching historical patient is found. Finally, the synthetic duration is sampled at random from the resulting pool of historical observations.

\section{Fairness-Aware Optimization Model}
\label{sec:FAIR:opt}

In this section, we present the fairness-aware MILP model used to compute outpatient schedules. We first formulate the outpatient scheduling problem, then describe the proposed optimization model. Next, we introduce an $\varepsilon$-constraint reformulation together with a dominance reduction strategy that limits the number of $\varepsilon$ values that must be evaluated. Finally, we present a baseline model without fairness considerations to assess the effectiveness of the proposed approach.

\subsection{Problem Definition}

We formulate the outpatient scheduling problem as the joint assignment of patients to examination days, CT scanners, and radiologists over a one-week planning horizon. Specifically, the model determines (i) the examination day assigned to each patient, (ii) the CT scanner on which the examination is performed, and (iii) the radiologist responsible for reporting the examination. Although CT scanners are assumed to be interchangeable, radiologists must be assigned according to their clinical specialty and patient--radiologist compatibility constraints. These compatibility requirements may naturally lead to unbalanced workloads; therefore, the proposed model seeks to distribute reporting activities as fairly as possible while satisfying all resource capacity constraints.

More generally, the proposed framework allows patient-specific scheduling preferences to be represented through a penalty associated with each feasible appointment slot (day and shift). Such penalties may reflect clinical recommendations, organizational requirements, or patient availability. The optimization model seeks to minimize the overall scheduling penalty while balancing radiologist workloads. In this paper, for simplicity, we assume that each patient is associated with a single clinically preferred examination day and define the penalty as the absolute deviation, measured in days, from this target date. For simplicity, and because no meaningful differences in urgency were observed among follow-up examinations in the case study, we do not distinguish between different urgency levels. Nevertheless, the proposed framework can be readily extended by incorporating urgency-dependent penalty functions.

The optimization focuses exclusively on outpatient scheduling, since emergency patients and inpatients are not scheduled in advance. Moreover, only the morning shift is considered because, in the case study, outpatient examinations are performed exclusively during this period, which is also the least congested. Consequently, shifts and days are equivalent in the proposed framework and will be referred to simply as days throughout the remainder of the paper.

The available capacity of each CT scanner is assumed to coincide with the duration of the morning shift. Conversely, the reporting capacity of each radiologist may exceed the nominal shift duration, reflecting current clinical practice. While reports for emergency patients and inpatients must be completed during the same working shift, even if overtime is required, outpatient reports can be finalized later when reporting capacity becomes available. Consequently, the proposed model balances the reporting workload assigned within each shift while allowing part of the outpatient reporting activity to be completed outside regular working hours.

Finally, the available capacity of both CT scanners and radiologists is reduced by the expected workload generated by emergency patients and inpatients. This adjustment accounts for both examination and reporting activities, ensuring that only the residual capacity is allocated to outpatient examinations.

To capture the variability in examination and reporting durations across patients, the optimization model incorporates patient--radiologist-specific predictions generated by the ML models presented in Section \ref{sec:FAIR:prediction}.

\subsection{Fairness-aware model}

\begin{table}[htbp]
\centering
\renewcommand{\arraystretch}{1.15}
\begin{tabular}{lp{15cm}}
\hline
\multicolumn{2}{l}{\textbf{Sets}}\\
\hline
$\mathcal{P}$ & Set of patients \\
$\mathcal{J}$ & Set of specialties \\
$\mathcal{R}$ & Set of radiologists \\
$\mathcal{R}_j$ & Set of radiologists of specialty $j$\\
$\mathcal{C}$ & Set of CT scanners \\
$\mathcal{D}$ & Set of days \\
\hline
\multicolumn{2}{l}{\textbf{Parameters}}\\
\hline
$\zeta_{pd}$ & Expected examination duration for patient $p$ in day $d$ (the value depends on the radiologist assigned to their specialty on that day)\\
$\nu_{pd}$ & Expected reporting duration for patient $p$ in day $d$ (the value depends on the radiologist assigned to their specialty on that day)\\
$\sigma_d$ & Expected total examination time of inpatient and emergency patients in day $d$\\
$\gamma_{dj}$ & Expected total reporting time associated with inpatient and emergency patients of specialty $j$ in day $d$\\
$\alpha_{pr}$  & 1 if the report of patient $p$ can be assigned to radiologist $r$, 0 otherwise \\
$\lambda_{pd}$ & 1 if day $d$ is the preferred day for patient $p$'s examination \\
$j_r$ & Specialty of radiologist $r$ \\
$\Xi$ & Shift duration\\
$\Xi^{\text{ref}}$ & Time available to each radiologist for completing reports in each shift\\
\hline
\multicolumn{2}{l}{\textbf{Decision Variables}}\\
\hline
$X_{prd}$ & 1 if patient $p$ is assigned to radiologist $r$ in day $d$, 0 otherwise \\
$Y_{pcd}$ & 1 if patient $p$ is assigned to CT scan $c$ in day $d$, 0 otherwise \\
$Z_p$ & Integer value indicating the number of days between the target day and the assigned examination day of patient $p$\\
$W_{rd}$ & Total workload of radiologist $r$ in day $d$\\

\hline
\end{tabular}
\caption{Notation of sets, parameters, and decision variables in the mathematical 
         formulation.}
\label{tab:FAIR:sets_parameters}
\end{table}

In the proposed model, fairness is considered from both the radiologists' and the patients' perspectives. We formulate a multi-objective optimization problem that simultaneously minimizes the reporting workload imbalance among radiologists and the deviation between patients' clinically preferred and scheduled examination days. For both objectives, fairness is enforced through a lexicographic min--max criterion \cite{Bertsimas2011}, a well-established approach in the fairness optimization literature. Such a criterion generalizes the principles of Rawlsian justice and the Kalai--Smorodinsky bargaining solution to settings involving multiple decision-makers. Unlike aggregate measures of fairness, it successively minimizes the largest individual disadvantage, then the second largest, and so forth, preventing improvements in average performance from masking highly unfavorable outcomes for individual radiologists or patients. Sets, parameters, and decision variables of the mathematical formulation are defined in Table~\ref{tab:FAIR:sets_parameters}. 

The model constraints are formulated as follows
\begin{align}
     \label{eq:FAIR:one_ct}
    \sum_{c\in\mathcal{C}}\sum_{d\in\mathcal{D}} Y_{pcd} & = 1, \quad &&\forall p\in\mathcal{P},\\
    \label{eq:FAIR:shift_duration}
    \sum_{p\in\mathcal{P}} \zeta_{pd} Y_{pcd}&\le \Xi-\bigg\lceil\frac{\sigma_d}{|\mathcal{C}|}\bigg\rceil, &&\forall c\in\mathcal{C},d\in\mathcal{D},\\
    \label{eq:FAIR:one_day}
    \sum_{r\in\mathcal{R}}\sum_{d\in\mathcal{D}} X_{prd} &= 1,  &&\forall p\in\mathcal{P},\\
     \label{eq:FAIR:rep_duration}
     \sum_{p\in\mathcal{P}} \nu_{pd} X_{prd}&\le \Xi^{ref} - \bigg\lceil\frac{\gamma_{dj_r}}{|\mathcal{R}_{j_r}|}\bigg\rceil, &&\forall r\in\mathcal{R},d\in\mathcal{D},\\
     \label{eq:FAIR:rad_equal_day}
     \sum_{r\in\mathcal{R}} X_{prd} &= \sum_{c\in\mathcal{C}}Y_{pcd}, &&\forall p\in\mathcal{P}, d\in\mathcal{D},\\
     \label{eq:FAIR:compatibility}
     X_{prd}& \le \alpha_{pr}, &&\forall p\in\mathcal{P}, r\in\mathcal{R},d\in\mathcal{D},\\
     \label{eq:FAIR:delay_days}
      Z_p &= \left| \sum_{d\in\mathcal{D}} d\,\lambda_{pd} - \sum_{d\in\mathcal{D}} \sum_{c\in\mathcal{C}} d\,Y_{pcd} \right|, &&\forall p\in\mathcal{P},\\
      \label{eq:FAIR:workload}
      W_{rd} &\ge \sum_{p\in \mathcal{P}}\nu_{pd} X_{prd} + \bigg\lceil\frac{\gamma_{dj_r}}{|\mathcal{R}_{j_r}|}\bigg\rceil, &&\forall r\in \mathcal{R}, \forall d\in \mathcal{D},\\
      \label{eq:FAIR:domains1}
    X_{prd}, Y_{pcd} &\in \{0,1\}, &&\forall p\in\mathcal{P}, r\in\mathcal{R}, c\in\mathcal{C}, d\in\mathcal{D},\\
    Z_p&\in\mathbb{N}_0, && \forall p\in\mathcal{P}, \\
    \label{eq:FAIR:domains4}
    Z_p &\le |\mathcal{D}|-1, && \forall p\in\mathcal{P},\\
    \label{eq:FAIR:domains3}
    W_{rd}& \ge 0, &&\forall r\in\mathcal{R}, d\in\mathcal{D},\\
    \label{eq:FAIR:domains2}
W_{rd} &\le  \Xi^{ref}, &&\forall r\in\mathcal{R}, d\in\mathcal{D}.
\end{align}

Constraints~\eqref{eq:FAIR:one_ct} ensure that each patient is assigned to exactly one CT scanner and one examination day. Constraints~\eqref{eq:FAIR:shift_duration} impose that the expected total examination time of the patients assigned to a CT scanner on a given day does not exceed the available shift duration. Constraints~\eqref{eq:FAIR:one_day} guarantee that each patient is assigned to exactly one radiologist on the same examination day. Constraints~\eqref{eq:FAIR:rep_duration} ensure that the expected total reporting time of the patients assigned to a radiologist on a given day does not exceed the maximum radiologist workload, where $j_r$ denotes the specialty of radiologist $r$, and $\mathcal{R}_j$ is the set of radiologists belonging to specialty $j$. 
Constraints~\eqref{eq:FAIR:rad_equal_day} impose that the day on which a patient undergoes the CT examination coincides with the day of the radiologist assigned to interpret the examination. Constraints~\eqref{eq:FAIR:compatibility} guarantee patient--radiologist compatibility requirements. Constraints~\eqref{eq:FAIR:delay_days} define the number of days between the patient's preferred examination day and the actual scheduled examination day. After linearization, they become:
\begin{align}
Z_p &\ge \sum_{d\in\mathcal{D}} d\lambda_{pd} - \sum_{d\in\mathcal{D}} \sum_{c\in\mathcal{C}} dY_{pcd}, && \quad \forall p\in\mathcal{P},\\
Z_p &\ge \sum_{d\in\mathcal{D}} \sum_{c\in\mathcal{C}} dY_{pcd} - \sum_{d\in\mathcal{D}} d\lambda_{pd}, && \quad \forall p\in\mathcal{P}.
\end{align}
Constraints~\eqref{eq:FAIR:workload} define the workload of radiologist $r$ on day $d$, while constraints~\eqref{eq:FAIR:domains1}--\eqref{eq:FAIR:domains2} define the domains of decision variables.

The radiologists' objective is to minimize the maximum workload in each day, while the patients' objective is to minimize the maximum wait
    \begin{equation}
        \min \max_{r\in\mathcal{R},\, d\in\mathcal{D}} W_{rd}, \quad \quad \min \max_{p\in\mathcal{P}} Z_{p},
    \end{equation}
    which can be reformulated as
    \begin{equation}
    \begin{aligned}
        \min \quad & t \\
        \text{s.t.} \quad & t \ge W_{rd} \quad \forall r\in\mathcal{R},\, d\in\mathcal{D},
    \end{aligned}
    \quad \quad 
     \begin{aligned}
        \min \quad & s \\
        \text{s.t.} \quad & s \ge Z_{p} \quad \forall p\in\mathcal{P}.
    \end{aligned}
    \end{equation}

Following the criterion formalized by \cite{Bertsimas2011}, both the radiologists' and the patients' objectives are optimized lexicographically. The procedure first minimizes the maximum workload (or waiting time). The optimal value is then fixed, and the second-largest workload (or waiting time) is minimized. The process continues iteratively until all workloads (or waiting times) have been determined.

\subsection{Generalization of the \texorpdfstring{$\varepsilon$}{epsilon}-constraint method to the lexicographic min--max}

We adopt a generalized $\varepsilon$-constraint approach, where the patient waiting time distribution is controlled via a vector of thresholds. The primary objective is to minimize the
maximum reporting workload among all radiologists and shifts
\begin{equation}
\min \max_{r\in\mathcal{R},\, d\in\mathcal{D}} W_{rd}.
\end{equation}

To control patient waiting times, we introduce a constraint based on a representation of waiting-time distributions. Let $\varepsilon \in \Omega$, where $\Omega$ denotes the set of all feasible
waiting-time allocations over $\mathcal{D}$ days, defined as
vectors $\varepsilon = (\varepsilon_d)_{d\in\mathcal{D}}$ such that $\varepsilon_d \in \mathbb{N}_0$ and
\[
\sum_{d\in\mathcal{D}} \varepsilon_d = |\mathcal{P}|.
\]
Each component $\varepsilon_d$ represents the maximum number of patients allowed to wait exactly $d$ days.

In addition, we introduce the following auxiliary variables:
\begin{align*}
K_{pd} &\in \{0,1\}
&& \text{equal to 1 if patient $p$ is assigned a waiting time of $d$ days},\\
H_d &\in \mathbb{N}_0
&& \text{number of patients assigned a waiting time of $d$ days}.
\end{align*}
These variables are defined through the following constraints
\begin{align}
\label{eq:FAIR:def_k}
Z_p &= \sum_{d\in\mathcal{D}} dK_{pd},
&& \forall p\in\mathcal{P},\\
\sum_{d\in\mathcal{D}} K_{pd} &= 1,
&& \forall p\in\mathcal{P},\\
\label{eq:FAIR:def_h}
H_d &= \sum_{p\in\mathcal{P}} K_{pd},
&& \forall d\in\mathcal{D}.
\end{align}

Finally, for a given $\varepsilon \in \Omega$, the constraints
\begin{equation}
\label{eq:FAIR:bound_wt}
H_d \le \varepsilon_d,
\qquad \forall d\in\mathcal{D}
\end{equation}
ensures that the waiting-times of patients does not exceed $\varepsilon$.

For each fixed $\varepsilon$, lexicographic min--max fairness is achieved by iteratively minimizing the radiologists' workloads. First, the maximum workload is minimized. Once its optimal value is obtained, it is fixed via an additional constraint, and the second-largest workload is minimized. This procedure is repeated by successively fixing the previously optimized workload levels and minimizing the next-largest one. The complete algorithm is reported in Algorithm~\ref{alg:FAIR:minmax}.

\begin{algorithm}[htb]
\label{algo:FAIR:minmax}
\caption{Lexicographic Min--Max Fairness}
\label{alg:FAIR:minmax}
\begin{algorithmic}[1]
\State \textbf{Input:} Instance characteristics ($\mathcal{P},\mathcal{J},\mathcal{R},\mathcal{C}, \mathcal{D}, \zeta, \nu, \alpha, \lambda,\Xi, \Xi^{ref}$), $\varepsilon$
\State \textbf{Output:} Assignment of patients to CT scanners and radiologists that satisfies the lexicographically minimal workload $(X^*, Y^*, W^*)$, and corresponding workloads $b_{rd}$ for all $r\in\mathcal{R}, d\in\mathcal{D}$

\State $\mathcal{F}_{\text{free}} \leftarrow \mathcal{R}\times\mathcal{D}$ \Comment{set of (radiologist, day) pairs not yet fixed}
\State $\mathcal{F}_{\text{fixed}} \leftarrow \emptyset$ \Comment{set of (radiologist, day) pairs fixed}
\State $b_{rd} \leftarrow +\infty \quad \forall (r,d)\in\mathcal{R}\times\mathcal{D}$ \Comment{fixed workload upper bounds}

\While{$\mathcal{F}_{\text{free}} \neq \emptyset$}
    \State \textbf{Solve MILP:}
    \State  $\min t$ s.t. \eqref{eq:FAIR:one_ct}--\eqref{eq:FAIR:domains2}, \eqref{eq:FAIR:def_k}--\eqref{eq:FAIR:bound_wt}, $W_{rd} \le t \, \forall (r,d)\in\mathcal{F}_{\text{free}}$, $W_{rd} \le b_{rd} \, \forall (r,d)\in\mathcal{F}_{\text{fixed}}$
    \State Let $t^*$ be the optimal value, and let $(X^*, Y^*, W^*)$ be the optimal solution, \Statex \quad \quad where $W^*_{rd}$ denotes the workload of radiologist $r$ in day $d$

    \State $\text{Saturated} \leftarrow \emptyset$
    \For{each $(r,d) \in \mathcal{F}_{\text{free}}$}
        \If{$W^*_{rd} = t^*$}
            \State $\text{Saturated} \leftarrow \text{Saturated} \cup \{(r,d)\}$
        \EndIf
    \EndFor


    \For{each $(r,d) \in \text{Saturated}$}
        \State $b_{rd} \leftarrow t^*$
        \State $\mathcal{F}_{\text{fixed}} \leftarrow \mathcal{F}_{\text{fixed}} \cup \{(r,d)\}$
        \State $\mathcal{F}_{\text{free}} \leftarrow \mathcal{F}_{\text{free}} \setminus \{(r,d)\}$
    \EndFor
\EndWhile

\State \textbf{return} $(X^*, Y^*, W^*)$ and $b_{rd} \quad \forall (r,d)\in\mathcal{R}\times\mathcal{D}$
\end{algorithmic}
\end{algorithm}

\subsection{Cardinality and dominance reduction}

The set $\Omega$ is in one-to-one correspondence with the set of distributions of
$|\mathcal{P}|$ indistinguishable items into $|\mathcal{D}|$ categories.
By the classical \emph{Stars and Bars} theorem~\cite{stanley2011ec1},
its cardinality is
\begin{equation}
|\Omega|
= \binom{|\mathcal{D}| + |\mathcal{P}| - 1}{|\mathcal{P}|}
= \frac{(|\mathcal{D}| + |\mathcal{P}| - 1)!}{|\mathcal{P}|!\, (|\mathcal{D}| - 1)!}.
\end{equation}

For instance, for $|\mathcal{P}|=15$ patients and $|\mathcal{D}|=5$ possible
waiting days, it holds $|\Omega|=3,876$. Exhaustively solving the MILP for
each $\varepsilon \in \Omega$ with a brute-force approach is therefore computationally expensive.
To reduce the number of evaluations, we exploit a monotonicity property of
the feasible region. Specifically, for $\varepsilon, \hat{\varepsilon} \in \Omega$,
we define the partial order
\[
\varepsilon \preceq \hat{\varepsilon}
\quad \Longleftrightarrow \quad
\varepsilon_d \le \hat{\varepsilon}_d \ \ \forall d \in \mathcal{D}.
\]
In this case, $\hat{\varepsilon}$ is said to dominate $\varepsilon$.

Since the feasible region is monotone in $\varepsilon$, if the model is
infeasible for some $\hat{\varepsilon}$, then it is also infeasible for all
$\varepsilon \preceq \hat{\varepsilon}$. Therefore, dominated vectors of an
unfeasible $\varepsilon$ need not be evaluated.

Moreover, if two vectors $\varepsilon$ and $\hat{\varepsilon}$ yield the same
optimal objective value (i.e., they have the same workloads for each couple radiologist-shift) and $\varepsilon \preceq \hat{\varepsilon}$, then all
vectors $\varepsilon'$ such that
\[
\varepsilon \preceq \varepsilon' \preceq \hat{\varepsilon}
\]
do not generate improved Pareto-optimal solutions and can be safely discarded
from the enumeration. 
This dominance-based pruning significantly reduces the number of MILP
instances that must be solved while preserving the Pareto frontier.

\subsection{Fairness-unaware baseline model}

As a benchmark, we consider an aggregate-efficiency baseline model that minimizes the total excess reporting workload and the total deviation between scheduled and preferred examination days. The objective function is given by
\begin{equation}
\min \sum_{d\in\mathcal{D}}\sum_{r\in\mathcal{R}}(W_{rd}-\Xi)+\sum_{p\in\mathcal{P}}Z_p.
\end{equation}

The model is subject to constraints \eqref{eq:FAIR:one_ct}--\eqref{eq:FAIR:domains2}.
As in the fairness-aware formulation, we adopt an $\varepsilon$-constraint approach. Specifically, we consider the objective function
\begin{equation}
\min \sum_{d\in\mathcal{D}}\sum_{r\in\mathcal{R}}(W_{rd}-\Xi)
\end{equation}
and introduce the additional constraints
\begin{equation}
\sum_{p\in\mathcal{P}}Z_p \le \tilde{\varepsilon}, \quad \forall \tilde{\varepsilon}\in\{0,\dots,|\mathcal{D}|\cdot|\mathcal{I}|\}.
\end{equation}

In this setting, the number of MILPs to be solved is limited, and the corresponding models can be totally ordered according to the value of $\tilde{\varepsilon}$. Therefore, we solve the models starting from the largest value of $\tilde{\varepsilon}$ and proceed in decreasing order. As soon as an instance becomes infeasible, we terminate the procedure, since all models associated with smaller values of $\tilde{\varepsilon}$ are necessarily infeasible as well.

\section{Data Description, Results and Discussion}
\label{sec:FAIR:results}

This section outlines the data infrastructure and pre-processing procedures adopted in the study, and the results obtained from the predictive models and the optimization framework, along with managerial insight.

All experiments were conducted in \texttt{Python 3.13}. ML and deep learning models were implemented using the \texttt{Scikit-learn} library \cite{Pedregosa}.

\subsection{Problem Description and Scheduling Framework}

The study was conducted at the ED of Fondazione IRCCS Policlinico San Matteo, a tertiary trauma center in Northern Italy serving a catchment population of approximately 600,000 to 1.2 million inhabitants. The hospital represents one of the busiest facilities within the Italian emergency care network. The radiology department operates two CT scanners during the daytime shifts (08:00--14:00 and 14:00--20:00) and one CT scanner overnight (20:00--08:00). Although the scanners differ slightly in their technical characteristics, they are interchangeable for almost all examinations. During each daytime shift, two radiologists are on duty, one specializing in body imaging (shortly \textit{body}) and the other in neuroradiology (shortly \textit{neuro}).

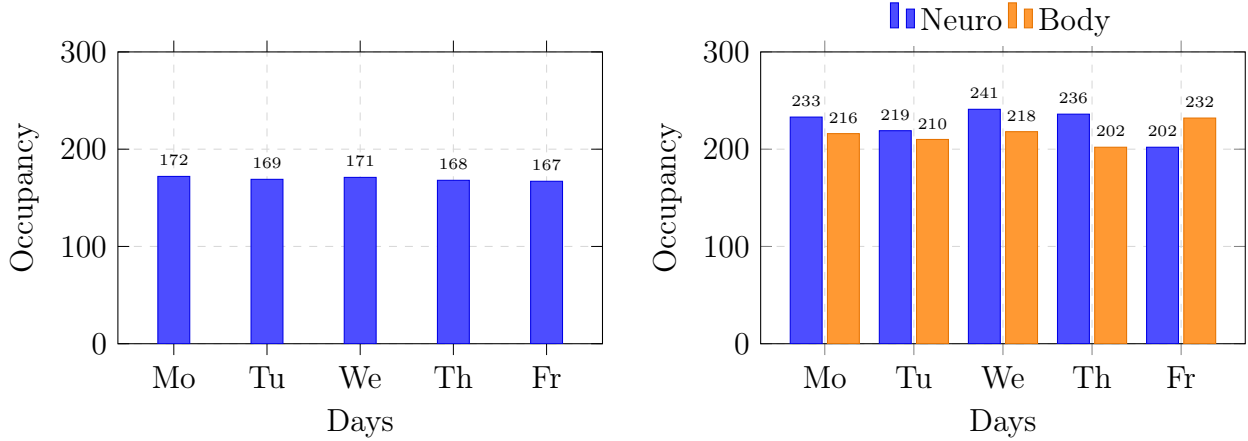
\begin{figure}[t]
\centering

\begin{tikzpicture}

\begin{axis}[
    at={(-3cm,0)},
    width=0.47\textwidth,
    height=0.32\textwidth,
    ybar,
    bar width=12pt,
    ymin=0,
    ymax=300,
    xlabel={Days},
    ylabel={Occupancy},
    symbolic x coords={Mo,Tu,We,Th,Fr},
    xtick=data,
    grid=major,
    grid style={dashed,gray!30},
    axis line style={black},
    tick style={black},
    enlarge x limits=0.15,
    nodes near coords,
    every node near coord/.append style={font=\tiny},
]

\addplot[
    fill=blue!70,
    draw=blue!90!black
] coordinates {
    (Mo,172)
    (Tu,169)
    (We,171)
    (Th,168)
    (Fr,167)
};

\end{axis}

\begin{axis}[
    at={(5.5cm,0)},
    anchor=south west,
    width=0.47\textwidth,
    height=0.32\textwidth,
    ybar,
    bar width=12pt,
    ymin=0,
    ymax=300,
    xlabel={Days},
    ylabel={Occupancy},
    symbolic x coords={Mo,Tu,We,Th,Fr},
    xtick=data,
    grid=major,
    grid style={dashed,gray!30},
    enlarge x limits=0.18,
    legend style={
        at={(0.5,1.02)},
        anchor=south,
        legend columns=2,
        draw=none
    },
    nodes near coords,
    every node near coord/.append style={font=\tiny},
]

\addplot[
    fill=blue!70,
    draw=blue!90!black
] coordinates {
    (Mo,233)
    (Tu,219)
    (We,241)
    (Th,236)
    (Fr,202)
};

\addplot[
    fill=orange!80,
    draw=orange!90!black
] coordinates {
    (Mo,216)
    (Tu,210)
    (We,218)
    (Th,202)
    (Fr,232)
};

\legend{Neuro, Body}

\end{axis}

\end{tikzpicture}

\caption{Occupancy statistics. Left: mean occupancy of the CT scanners throughout the week. Right: mean occupancy of the two radiologists throughout the week.}
\label{fig:FAIR:occupancy}
\end{figure}

While outpatient examinations constitute the schedulable component of the CT workload, emergency and inpatient requests arrive unpredictably and frequently disrupt planned activities. The occupancy statistics of non-elective patients, both for CT scanners and radiologists, are shown in Figure~\ref{fig:FAIR:occupancy}. They correspond to the parameters $\sigma_d$ and $\gamma_{dj}$, respectively.

Table~\ref{tab:FAIR:patientDistributions} reports the distribution of examination specialties and urgency levels by patient category. As expected, emergency patients and inpatients are associated with substantially higher urgency levels than outpatients.

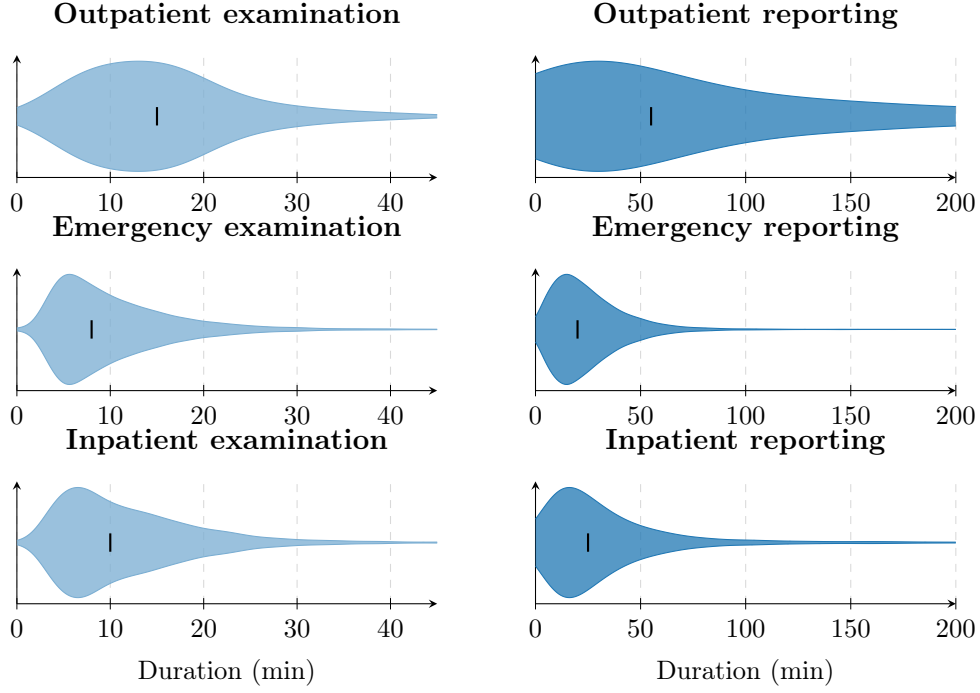
\begin{figure}[htbp]
\centering
\begin{tikzpicture}
\begin{groupplot}[
    group style={
        group size=2 by 3,
        horizontal sep=1.3cm,
        vertical sep=1.2cm,
    },
    width=0.42\textwidth,
    height=3.2cm,
    ymin=0, ymax=2,
    ytick=\empty,
    xlabel style={font=\footnotesize},
    axis lines=left,
    axis line style={black},
    tick style={black},
    tick label style={font=\footnotesize},
    title style={font=\small\bfseries},
    xmajorgrids,
    grid style={dashed,gray!30},
    every axis plot/.append style={no markers},
    clip=false,
]

\nextgroupplot[title={Outpatient examination}, xmin=0, xmax=45]
\addplot[name path=u1,draw=colOut!60,thin] table[x=x,y expr=1+0.9*\thisrow{d}] {figures/outexam.dat};
\addplot[name path=l1,draw=colOut!60,thin] table[x=x,y expr=1-0.9*\thisrow{d}] {figures/outexam.dat};
\addplot[colOut!60,fill=colOut!60,fill opacity=0.75] fill between[of=u1 and l1];
\addplot[black,thick] coordinates {(\medoutexam,0.85) (\medoutexam,1.15)};

\nextgroupplot[title={Outpatient reporting}, xmin=0, xmax=200]
\addplot[name path=u2,draw=colOut,thin] table[x=x,y expr=1+0.9*\thisrow{d}] {figures/outrep.dat};
\addplot[name path=l2,draw=colOut,thin] table[x=x,y expr=1-0.9*\thisrow{d}] {figures/outrep.dat};
\addplot[colOut,fill=colOut,fill opacity=0.75] fill between[of=u2 and l2];
\addplot[black,thick] coordinates {(\medoutrep,0.85) (\medoutrep,1.15)};

\nextgroupplot[title={Emergency examination}, xmin=0, xmax=45]
\addplot[name path=u3,draw=colOut!60,thin] table[x=x,y expr=1+0.9*\thisrow{d}] {figures/emgexam.dat};
\addplot[name path=l3,draw=colOut!60,thin] table[x=x,y expr=1-0.9*\thisrow{d}] {figures/emgexam.dat};
\addplot[colOut!60,fill=colOut!60,fill opacity=0.75] fill between[of=u3 and l3];
\addplot[black,thick] coordinates {(\medemgexam,0.85) (\medemgexam,1.15)};

\nextgroupplot[title={Emergency reporting}, xmin=0, xmax=200]
\addplot[name path=u4,draw=colOut,thin] table[x=x,y expr=1+0.9*\thisrow{d}] {figures/emgrep.dat};
\addplot[name path=l4,draw=colOut,thin] table[x=x,y expr=1-0.9*\thisrow{d}] {figures/emgrep.dat};
\addplot[colOut,fill=colOut,fill opacity=0.75] fill between[of=u4 and l4];
\addplot[black,thick] coordinates {(\medemgrep,0.85) (\medemgrep,1.15)};

\nextgroupplot[title={Inpatient examination}, xmin=0, xmax=45, xlabel={Duration (min)}]
\addplot[name path=u5,draw=colOut!60,thin] table[x=x,y expr=1+0.9*\thisrow{d}] {figures/inpexam.dat};
\addplot[name path=l5,draw=colOut!60,thin] table[x=x,y expr=1-0.9*\thisrow{d}] {figures/inpexam.dat};
\addplot[colOut!60,fill=colOut!60,fill opacity=0.75] fill between[of=u5 and l5];
\addplot[black,thick] coordinates {(\medinpexam,0.85) (\medinpexam,1.15)};

\nextgroupplot[title={Inpatient reporting}, xmin=0, xmax=200, xlabel={Duration (min)}]
\addplot[name path=u6,draw=colOut,thin] table[x=x,y expr=1+0.9*\thisrow{d}] {figures/inprep.dat};
\addplot[name path=l6,draw=colOut,thin] table[x=x,y expr=1-0.9*\thisrow{d}] {figures/inprep.dat};
\addplot[colOut,fill=colOut,fill opacity=0.75] fill between[of=u6 and l6];
\addplot[black,thick] coordinates {(\medinprep,0.85) (\medinprep,1.15)};

\end{groupplot}
\end{tikzpicture}
\caption{Distribution of examination and reporting durations for outpatient, emergency, and inpatient cases. Left column: examination durations. Right column: reporting durations. The black vertical segment marks the median.}
\label{fig:FAIR:duration_violin}
\end{figure}

Figure~\ref{fig:FAIR:duration_violin} illustrates the distributions of examination and reporting durations for the three patient categories. Emergency examinations exhibit relatively concentrated execution and reporting times, whereas inpatient examinations show greater variability. Outpatient examinations display more homogeneous execution times but a markedly right-skewed distribution of reporting durations. This behavior reflects the operational practice whereby outpatient reports are often completed during idle periods between urgent cases, introducing additional variability in reporting times.

The study dataset included 26,945 CT examinations, corresponding to 49,813 individual scans, performed at the ED of Fondazione IRCCS Policlinico San Matteo, Pavia, Italy, between June 2021 and December 2022.

Overall, the mean examination duration was approximately 11.0 minutes (median: 9.0 minutes; standard deviation (SD): 7.9 minutes), whereas the mean reporting duration was approximately 30.0 minutes (median: 20.0 minutes; SD: 38.3 minutes). These observations motivate the development of predictive models capable of estimating patient-specific examination and reporting durations.

\begin{table}[t]
\centering
\caption{Proportion of exam specialties and urgency codes depending on the type of patient.}
\label{tab:FAIR:patientDistributions}
\begin{tabular}{lccc}
\toprule
 & \textbf{Emergencies} & \textbf{Inpatients} & \textbf{Outpatients} \\
\midrule

Neuro & 58.0\% & 50.1\% & 35.9\% \\
Body & 42.0\% & 49.9\% & 64.1\% \\
\midrule

Non/Low Urgent & 39.9\% & 53.5\% & 99.0\% \\
Intermediate Urgent & 53.4\% & 46.5\% & 1.0\% \\
Highly Critical & 6.7\% & - & - \\

\bottomrule
\end{tabular}
\end{table}

\subsubsection{Data preparation}

Clinical and operational data for each CT examination were extracted from the Radiology Information System (RIS). The available variables included patient demographics (age and sex), clinical setting (emergency, inpatient, or outpatient), urgency code, examination type, number of examinations, referring department, and the reporting radiologist. An important feature of the proposed framework is that the reporting radiologist can be explicitly incorporated into the predictive models, as the radiologists assigned to each shift are known at scheduling time. This makes it possible to estimate both examination and reporting durations for each patient--radiologist pair, thereby generating predictions that are tailored to both patient and radiologist characteristics.

\begin{figure}
\resizebox{\textwidth}{!}{%
\begin{tikzpicture}[
  font=\large,
  slot/.style={draw, minimum width=1cm, minimum height=0.7cm, font=\large},
  busy/.style={slot, fill=#1!40},
  free/.style={slot, fill=gray!15},
  patient/.style={circle, draw, minimum size=0.55cm, font=\large\bfseries, inner sep=1pt},
  arrow/.style={-{Stealth[length=2mm]}, thick},
  label/.style={font=\large}
]

\def\N{16}
\foreach \i in {1,...,\N} {
  \pgfmathsetmacro\prev{\i-1}
  \node[free] (s\i) at (\prev*1.1, 0) {};
}

\node[label, left=6pt of s1] {Slot (5 min)};

\node[busy=green] at (s3) {};
\node[busy=green] at (s4) {};
\node[busy=red]   at (s5) {};
\node[busy=red]   at (s6) {};
\node[busy=red]   at (s7) {};
\node[busy=red]   at (s8) {};
\node[busy=green] at (s9)  {};
\node[busy=green] at (s10) {};
\node[busy=green] at (s11) {};
\foreach \i in {12,...,\N} { \node[free] at (s\i) {}; }

\foreach \i in {1,...,\N} {
  \pgfmathsetmacro\prev{\i-1}
  \node[font=\normalsize, below=2pt] at (\prev*1.1, -0.35) {\i};
}

\node[patient, fill=green!50, above=18pt of s3] (pg) {OUT};
\draw[arrow, green!60!black] (pg.south) -- (s3.north);

\node[patient, fill=red!50, above=18pt of s5] (pr) {EME};
\draw[arrow, red!70!black] (pr.south) -- (s5.north);
\draw[arrow, orange!80!black, dashed, thick]
  ($(s4.north east)+(0,0.05)$) -- 
  node[above, text=orange!80!black] {preempted}
  ($(s5.north west)+(5,0.05)$);

\node[slot, fill=green!40, below=28pt of s3] (lg) {};
\node[label, right=3pt of lg] (o1) {Outpatient};
\node[slot, fill=red!40, right=120pt of lg] (lr) {};
\node[label, right=3pt of lr] {Emergency};
\node[slot, fill=gray!15, right=120pt of lr] (lf) {};
\node[label, right=3pt of lf] {Idle};

\draw[thick, -{Stealth}] (-0.6, 0) -- (\N*1.1 - 0.1, 0);

\end{tikzpicture}
}
\caption{Visual representation of how reporting times were modeled.}
\label{fig:FAIR:reportingModel}
\end{figure}
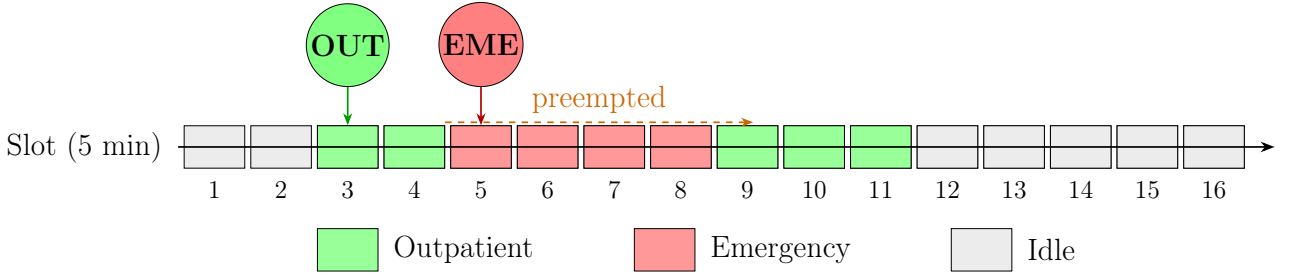

Two target variables were defined. The examination duration was computed as the time elapsed between patient entry into the CT scanner and the completion of image acquisition. The reporting duration was computed as the interval between the end of image acquisition and the validation of the final report in the RIS. Since reporting activities are frequently interrupted by other clinical duties, each 8-hour shift was discretized into five-minute intervals. Radiologists were assumed to be available for reporting during all intervals not occupied by examination-related activities, such as protocol definition, image quality assessment, or decisions regarding additional image acquisitions. Reporting followed a hierarchical priority policy: emergency examinations were always reported first, followed by inpatient and outpatient examinations, while examinations within each priority class were processed according to a first-come, first-served policy, as illustrated in Figure~\ref{fig:FAIR:reportingModel}: the reporting process for an outpatient case starts in slot 3. At slot 5, a reporting task for an emergency patient arrives, preempting the outpatient reporting process. The outpatient report remains suspended until the emergency report is completed in slot 8, after which the outpatient reporting process is resumed and completed.

To ensure data quality and reproducibility, several preprocessing steps were performed. Consecutive scans of the same patient conducted within 30 minutes were merged into a single examination, records with incomplete or inconsistent timestamps were removed, pediatric oncology examinations were excluded, and the dataset was restricted to examinations performed on the two ED CT scanners routinely used for emergency and inpatient imaging. The resulting dataset was subsequently used to train and test the predictive models and to provide the input data for the fairness-aware optimization model.
In particular, the last 20 weeks of data are used for validating the optimization model, while the remaining data are split into training and test sets, with 80\% allocated to training the ML models and 20\% allocated to testing.

\subsection{Predictive Model Results}

Predictive performance was evaluated using standard metrics such as the coefficient of determination (R$^2$), mean absolute error (MAE), and root mean squared error (RMSE). For each algorithm, the hyperparameters were optimized using \texttt{Optuna} with the objective of minimizing a normalized aggregate score combining these three performance metrics. Then, the selected hyperparameters of every ML model are used for the regret evaluation, as described in Section \ref{sec:FAIR:prediction}. 

Feature relevance was assessed using \texttt{SHAP} \cite{Lundberg2017}, which quantifies the contribution of each feature to an individual prediction by measuring how the prediction changes when the feature is included in the model.

Specifically, features were ranked according to their mean absolute \texttt{SHAP} value, normalized by the total importance, and their cumulative contribution was calculated. Only features accounting for the first 95\% of the cumulative \texttt{SHAP} importance were retained. 

\subsubsection{Feature analysis}
The most important features vary depending on the target variable (exam duration or reporting duration) and the prediction technique employed. 

Overall, for exam duration prediction, the most important features were specific examination types (\textit{Abdominal CT}, \textit{Intracranial CT Angiography}, \textit{Brain CT}, and \textit{Chest CT}), together with \textit{department}, \textit{number of examinations}, and \textit{patient origin}.

For reporting time prediction, the most influential features included several examination types (\textit{Brain CT}, \textit{Intracranial CT Angiography}, \textit{Abdominal CT}, and \textit{Pelvis and Sacroiliac Joints CT}), as well as \textit{number of examinations}, \textit{department}, \textit{patient origin}, \textit{exam specialty} (i.e., the specialty of the examination required by the patient), \textit{age}, \textit{scheduled outpatient}, \textit{emergency code}, \textit{referring physician}, and \textit{request hour}.

The feature importance analysis confirmed that, for exam duration prediction, the type of examination was the strongest predictor, particularly for more complex procedures such as abdominal and intracranial CT scans.

For reporting time prediction, the most relevant features were again the type of examination and the number of examinations, together with patient-related information and the referring physician.

These findings are consistent with clinical intuition and provide evidence supporting the face validity of the proposed prediction models.

\subsubsection{Predictive performance}
We retained some models trained on the full set of features and others trained on the features selected by \texttt{SHAP}, depending on their performance. The final results are reported in Table \ref{tab:FAIR:models_performance}.
CatBoost and MLP achieved the best performance for execution duration, with the highest R$^2$ (0.289) and lowest RMSE (6.54 minutes), while Random Forest has the lowest MAE (4.11 minutes). CatBoost performed best for reporting duration, with an  R$^2$ of 0.234, a RMSE of 30.89 minutes, and a MAE of 15.87. RF and XGBoost performed competitively, showing only minor differences from CatBoost and MLP and RF. In contrast, DTs consistently underperformed, confirming that methods used in previous studies are no longer optimal when more sophisticated techniques are available.
Importantly, the differences between top-performing models were small, suggesting that in practice several methods could be deployed successfully.

The regression models significantly outperformed the baseline historical averages for both execution and reporting times. For execution, the baseline exhibited a MAE of 5.47 minutes and a RMSE of 7.76 minutes. All tested models achieved lower errors. 
For reporting, the baseline had an MAE of 19.16 minutes and an RMSE of 35.31 minutes.

\begin{table}[tbp]
\centering
\begin{tabular}{lccccccc}
\toprule
& \multicolumn{3}{c}{\textbf{Exam Duration}}
& \multicolumn{3}{c}{\textbf{Reporting Duration}}
& \textbf{Regret} \\
\cmidrule(lr){2-4}
\cmidrule(lr){5-7}
\cmidrule(l){8-8}
\textbf{Model}
& $\mathbf{R^2}$
& \textbf{MAE}
& \textbf{RMSE}
& $\mathbf{R^2}$
& \textbf{MAE}
& \textbf{RMSE}
& \textbf{Mean} \\
\midrule
\textbf{BL}       & 0.000 & 5.47 & 7.76 & 0.000 & 19.16 & 35.31 & 0.22 \\
\textbf{DT}       & 0.277 & 4.14 & 6.60 & 0.169 & 16.83 & 32.18 & 0.19 \\
\textbf{RF}       & 0.286 & \textbf{4.11} & 6.56 & 0.229 & 16.17 & 30.99 & 0.18 \\
\textbf{XGBoost}  & 0.288 & 4.12 & 6.55 & 0.223 & 16.02 & 31.11 & \textbf{0.15} \\
\textbf{CatBoost} & \textbf{0.289} & 4.12 & \textbf{6.54} & \textbf{0.234} & \textbf{15.87} & \textbf{30.89} & 0.17 \\
\textbf{MLP}      & \textbf{0.289} & 4.13 & \textbf{6.54} & 0.217 & 16.45 & 31.23 & 0.16 \\
\bottomrule
\end{tabular}
\caption{Predictive performance of the ML algorithms in terms of conventional prediction metrics and optimization regret. MAE and RMSE are expressed in minutes. Regret is reported as the mean over the Pareto-dominating $\varepsilon$ values across the instances of the test set.}
\label{tab:FAIR:models_performance}
\end{table}

With respect to regret, measured as the average value of the regret function over the Pareto-dominating $\varepsilon$ values across the instances of the test set, XGBoost emerges as the best-performing method. The remaining predictive techniques exhibit similar performance and consistently outperform the baseline approach. These results confirm the findings in \cite{Daldossi2026}, namely that conventional prediction accuracy metrics fail to identify the most suitable ML technique when it is embedded within an optimization framework.

\subsection{Optimization Model Results}
Since the primary objective of the predictive component is to support the generation of optimized schedules, we selected XGBoost as the prediction technique. Although it does not achieve the best performance according to traditional predictive accuracy metrics, it outperforms the other approaches in terms of regret. This result indicates that, when integrated into the optimization framework, XGBoost provides more effective and robust decision-support capabilities.

We conducted a series of computational experiments to evaluate the performance of the proposed optimization model under different scenarios.

\begin{description}
\item[\textbf{Pareto trade-off analysis:}] We investigate the impact of different values of $\varepsilon \in \Omega$ on the resulting Pareto frontier. The objective of this analysis is to evaluate how radiologist workloads change as greater flexibility is introduced in patients' examination scheduling and to provide managerial insights into the trade-off between scheduling flexibility and workload balancing.

\item[\textbf{Scalability analysis:}] We assess the scalability of the proposed model under increasing patient demand and resource availability. This analysis aims to evaluate computational performance and solution quality as the problem size grows, representing larger radiology departments.

\item[\textbf{Patient case-mix analysis:}] We evaluate the model under different exam specialty distributions. The objective is to analyze how variations in patient--radiologist compatibility requirements affect the resulting schedules and workload balance.

\item[\textbf{Scheduling preference analysis:}] We investigate the model performance under different distributions of patients' preferred examination days. This analysis aims to understand how scheduling preferences influence the trade-off between patient waiting time and workload balancing.
\end{description}



\subsubsection{Dominance Reduction}

For the case study instances, which involve 5 scheduling days and 15 outpatients, a total of 3,976 MILPs would need to be solved, corresponding to one optimization problem for each feasible value of $\varepsilon$. 
The number of $\varepsilon$ values that can be discarded because they are dominated by an infeasible solution is negligible. Across the 10 test instances, this reduction amounts to only 3.4 MILPs on average.

In contrast, a large number of consecutive $\varepsilon$ values yield the same objective value, allowing a substantial reduction in the number of optimization problems that must be solved. On average, 3,211.7 MILPs can be skipped, corresponding to an 80.78\% reduction in computational effort.

These results demonstrate the effectiveness of the proposed dominance reduction strategy. Many different waiting-time distributions among patients lead to identical radiologist workload values, making it unnecessary to evaluate every feasible value of $\varepsilon$. Consequently, the proposed approach substantially reduces the computational burden while preserving the set of non-dominated solutions.

\subsubsection{Pareto trade-off analysis}
We investigate the impact of varying $\varepsilon \in \Omega$, i.e., the effect of modifying the maximum allowable deviation between a patient's preferred examination day and the scheduled day. In particular, we analyze the solutions generated by the proposed $\varepsilon$-constraint method for each instance.

The experiments were conducted on the real-world case study, involving 15 outpatients, 2 radiologists, 2 CT scanners, and a planning horizon of 5 working days. Each patient's preferred examination day was generated according to a discrete uniform distribution over the five available days.

Figure~\ref{fig:FAIR:heat_case_study} presents the heatmaps of the mean radiologist workloads across the instances of the test set obtained while minimizing the maximum radiologist workload for different values of $\varepsilon$. For each solution, radiologist--day workloads are first sorted in descending order and then averaged across the instances.

The smallest value of $\varepsilon$ for which the problem is feasible is $(12,3,0,0,0)$, indicating that 12 patients are scheduled on their preferred examination day, while the remaining three patients are shifted by one day. As $\varepsilon$ increases, different workload distributions emerge. Interestingly, these distributions do not follow the lexicographic ordering of $\varepsilon$: lexicographically larger values do not necessarily produce more balanced workloads, and vice versa. This behavior confirms that the relationship between patient scheduling flexibility and workload fairness is non-monotonic.

Overall, the workload of the busiest radiologist decreases from 576.6 minutes for $\varepsilon=(12,3,0,0,0)$ to 494.7 minutes, a value attained by several $\varepsilon$ configurations. Conversely, the workload of the least busy radiologist increases from 231.0 to 299.8 minutes, with the latter also achieved by multiple $\varepsilon$ values. The intermediate workloads exhibit smaller, yet still noticeable, variations across the different configurations.

These results show that allowing greater flexibility in patient scheduling can substantially improve workload balance among radiologists. At the same time, the existence of multiple $\varepsilon$ configurations yielding similar workload distributions suggests that different patient waiting-time patterns may provide equivalent operational performance. Consequently, the heatmap offers decision-makers a comprehensive view of the trade-off between adherence to patients' preferred examination days and workload balancing, enabling the selection of the scheduling policy that best matches the operational priorities of the radiology department.

\begin{figure}[tbp]
\centering
\includegraphics[width=0.85\linewidth]{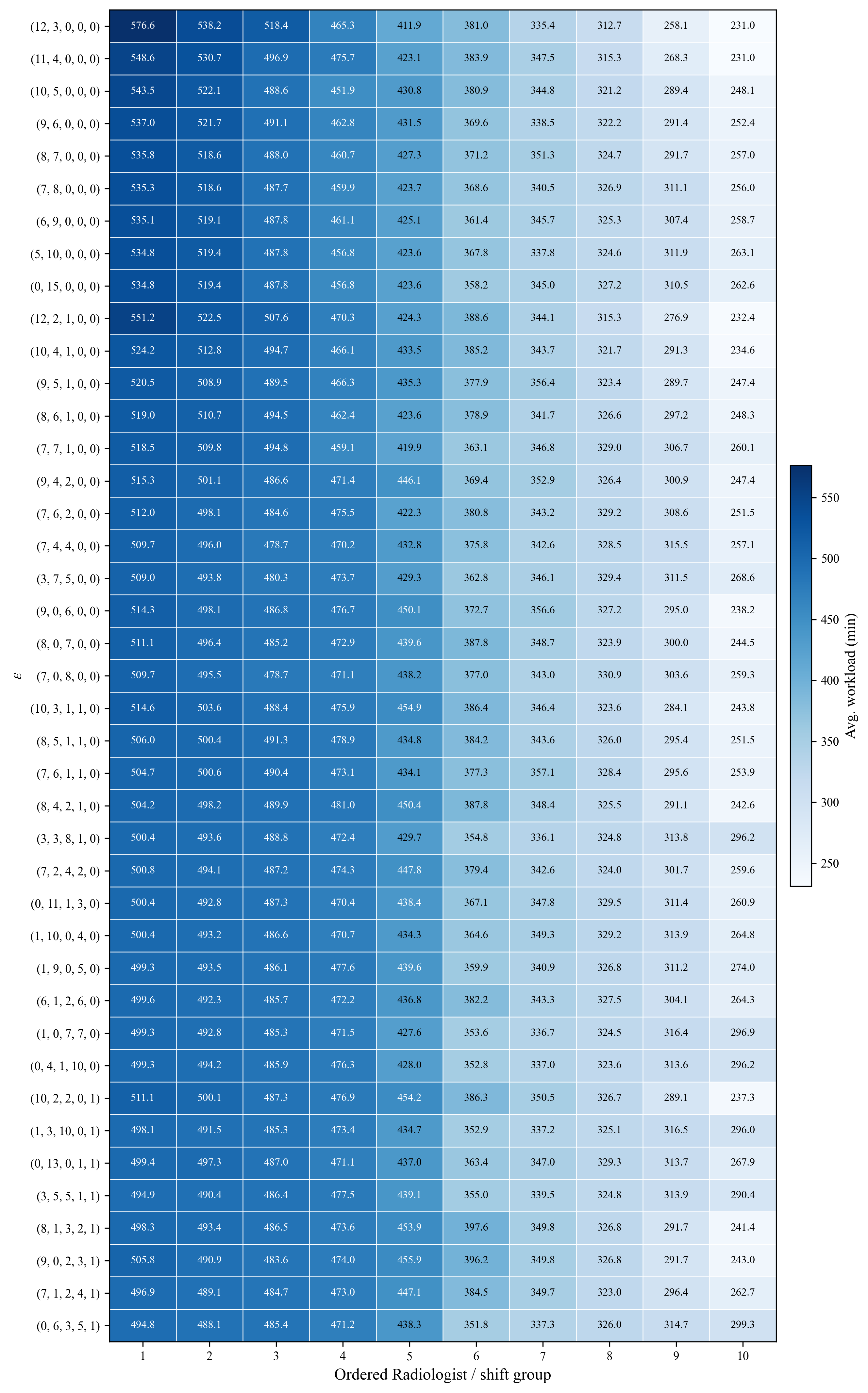}
\caption{Heatmap of the workload (in minutes) of radiologist-shift couples, depending on $\varepsilon$ (part 1).}
 \label{fig:FAIR:heat_case_study}
\end{figure}

\begin{figure}[tbp]
\ContinuedFloat
    \centering
    \includegraphics[width=0.85\linewidth]{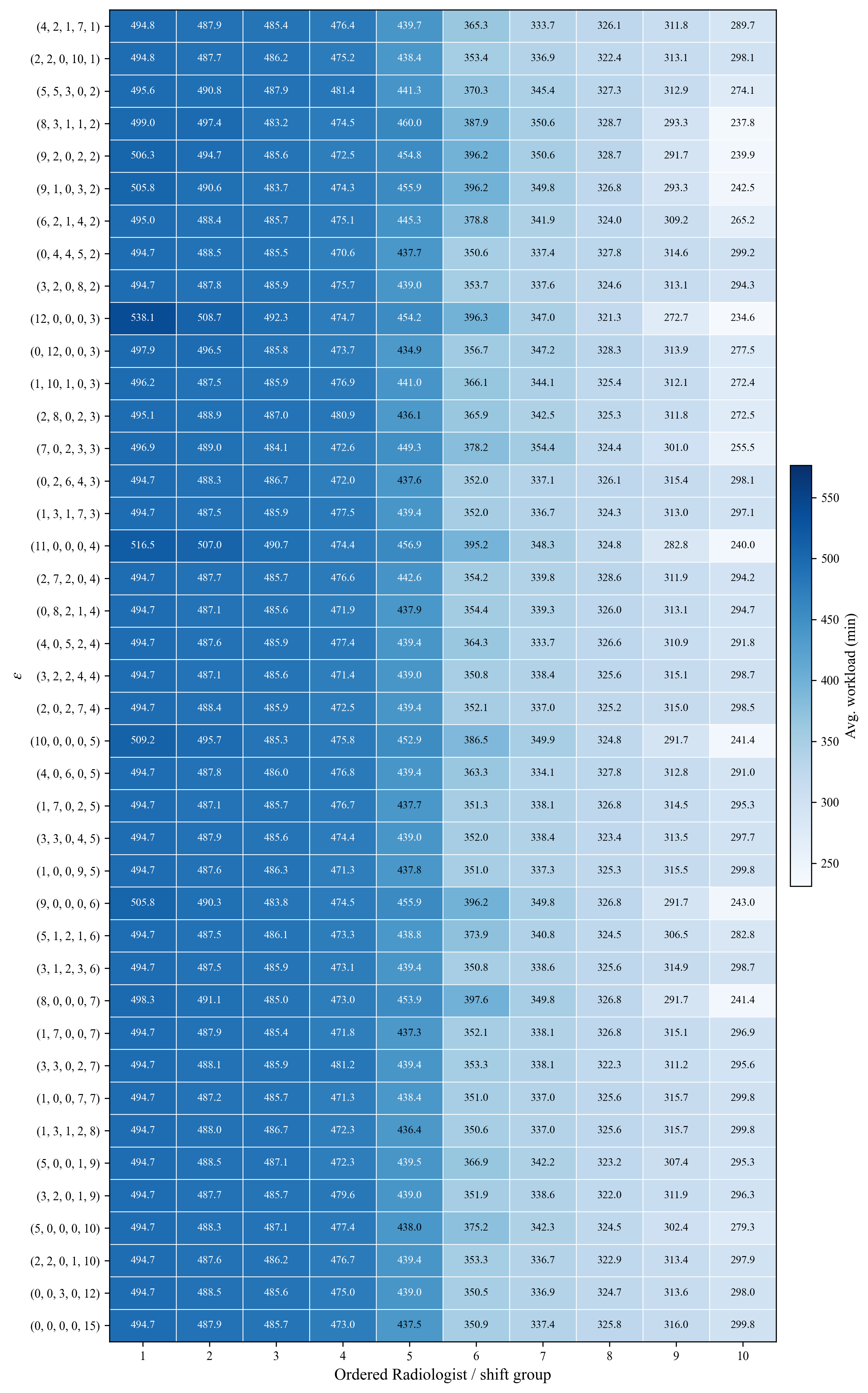}
    \caption{Heatmap of the workload (in minutes) of radiologist-shift couples, depending on $\varepsilon$ (part 2).}
\end{figure}

To gain further insight into the behavior of individual instances, we focus on a representative subset of $\varepsilon$ values. Specifically, we consider: (i) the smallest value for which all instances admit a feasible solution, namely $\varepsilon=(12,3,0,0,0)$; (ii) $\varepsilon=(0,15,0,0,0)$, which allows every patient to be scheduled up to one day after the preferred examination day; (iii) $\varepsilon=(0,0,15,0,0)$, which allows a delay of up to two days; (iv) $\varepsilon=(0,0,0,15,0)$, which allows a delay of up to three days; and (v) $\varepsilon=(0,0,0,0,15)$, which allows a delay of up to four days.

For each selected value of $\varepsilon$, Figure~\ref{fig:FAIR:spaghetti} reports boxplots of the ordered radiologist--day workloads (from the highest to the lowest), together with lines connecting the workload values obtained for each of the instances of the test set.

\begin{figure}[tbp]
\centering

\begin{subfigure}{0.75\textwidth}
    \centering
    \includegraphics[width=\linewidth]{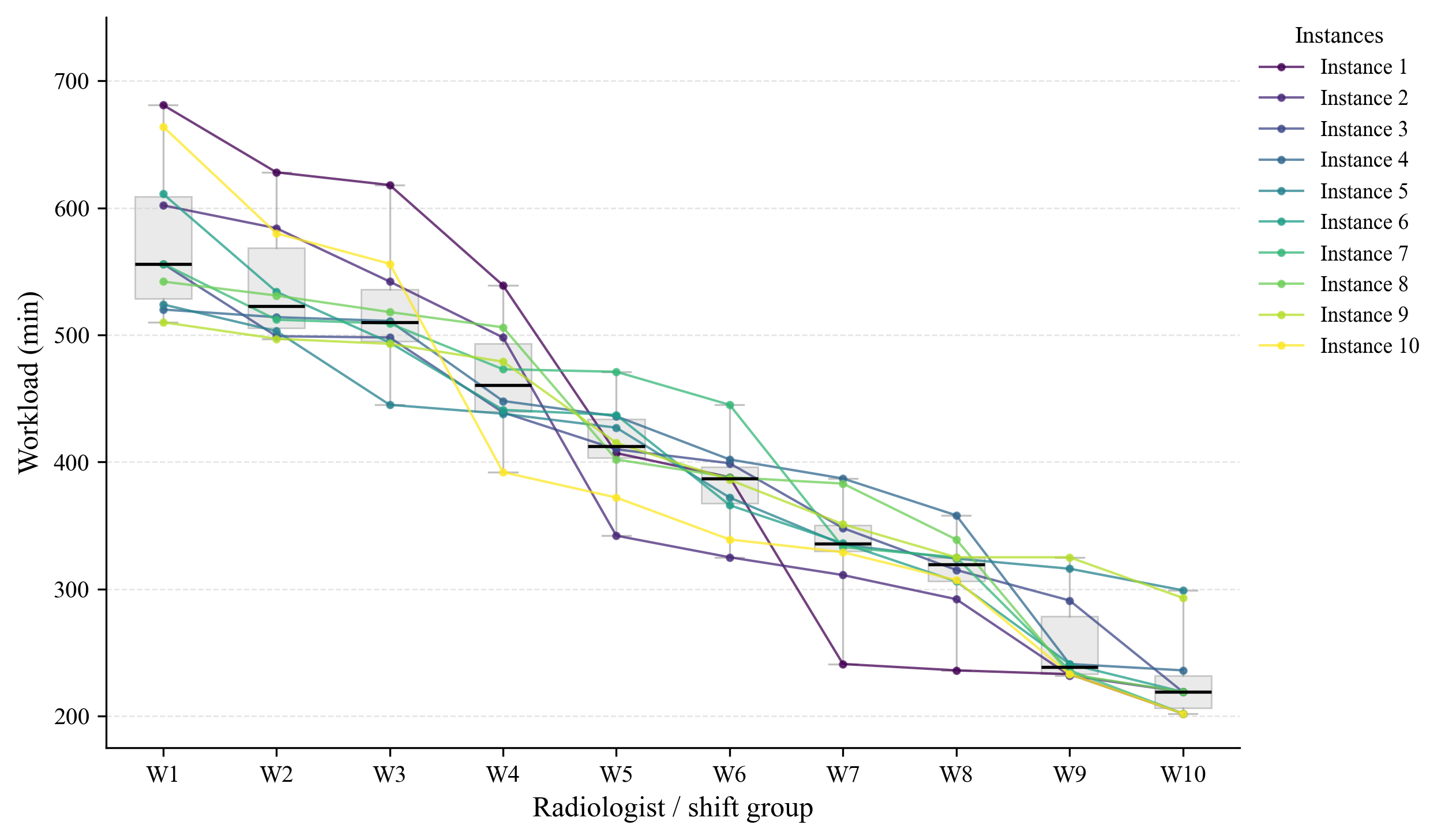}
    \caption{Workloads (in minutes) of the radiologists with $\varepsilon$ fixed as (i).}
    \label{fig:FAIR:inst1}
\end{subfigure}
\begin{subfigure}{0.75\textwidth}
    \centering
    \includegraphics[width=\linewidth]{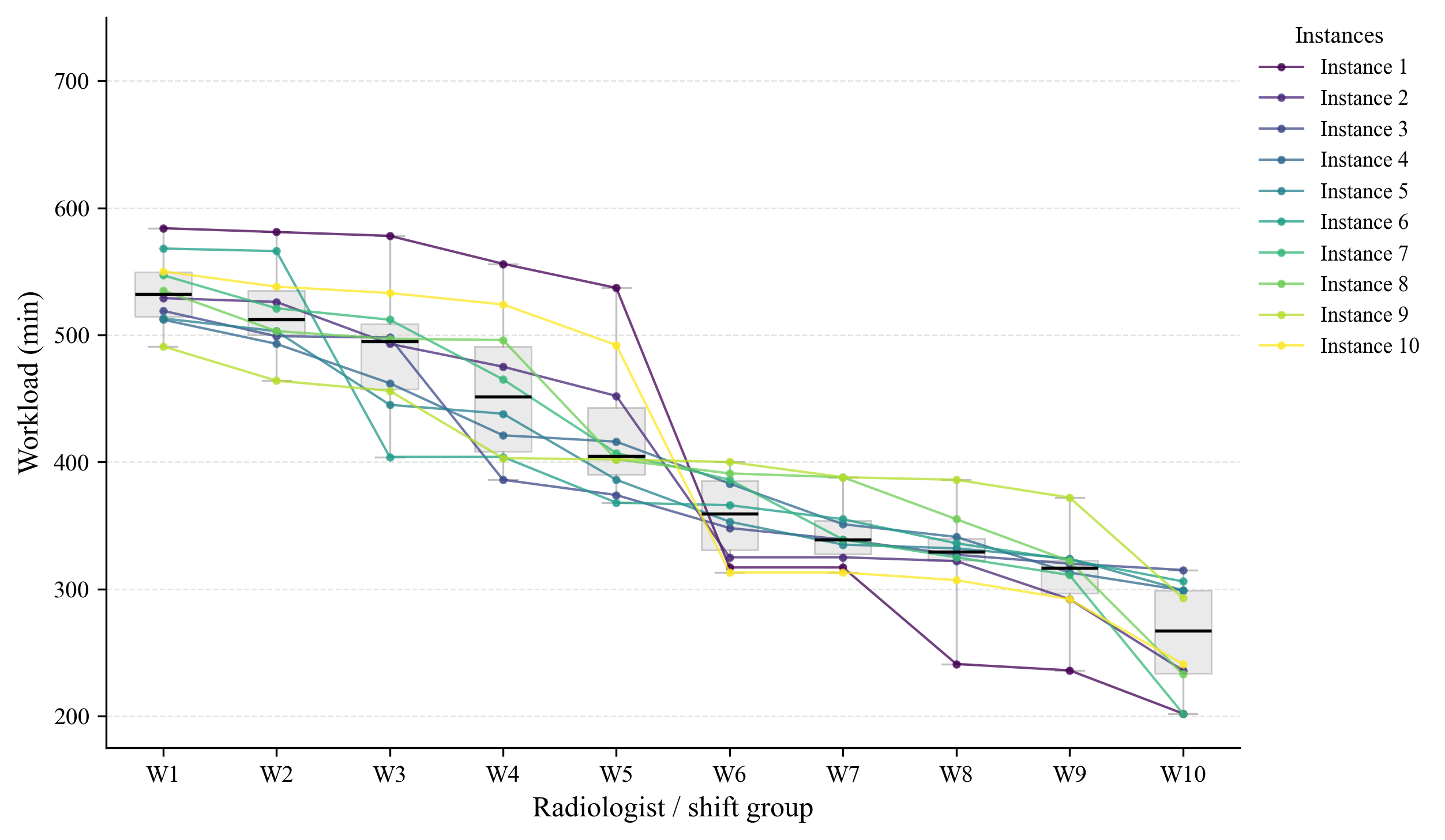}
    \caption{Workloads (in minutes) of the radiologists with $\varepsilon$ fixed as (ii).}
    \label{fig:FAIR:inst2}
\end{subfigure}
\begin{subfigure}{0.75\textwidth}
    \centering
    \includegraphics[width=\linewidth]{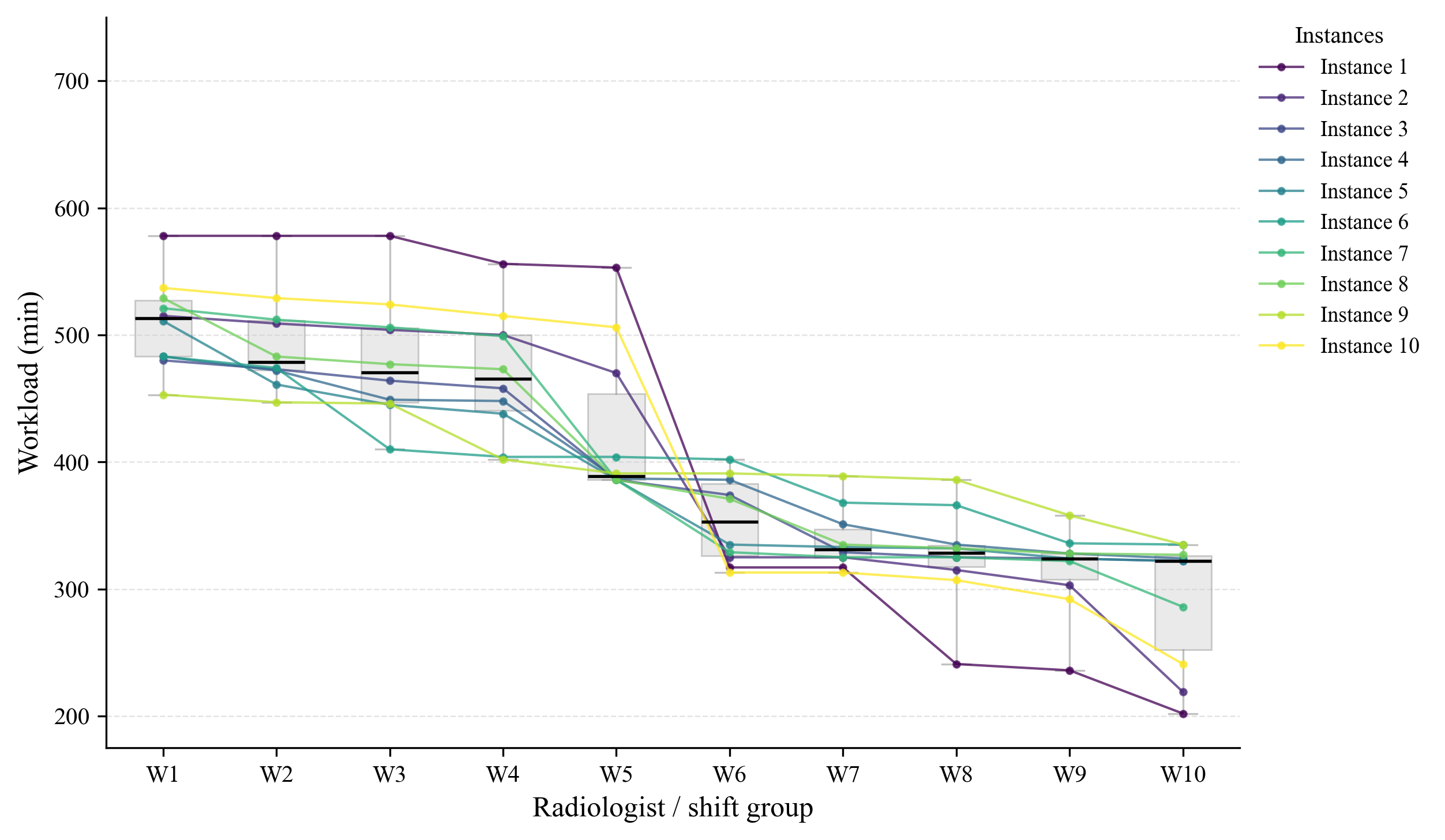}
    \caption{Workloads (in minutes) of the radiologists with $\varepsilon$ fixed as (iii).}
    \label{fig:FAIR:inst3}
\end{subfigure}
\end{figure}

\begin{figure}[tbp]
\ContinuedFloat
\centering

\begin{subfigure}{0.75\textwidth}
    \centering
    \includegraphics[width=\linewidth]{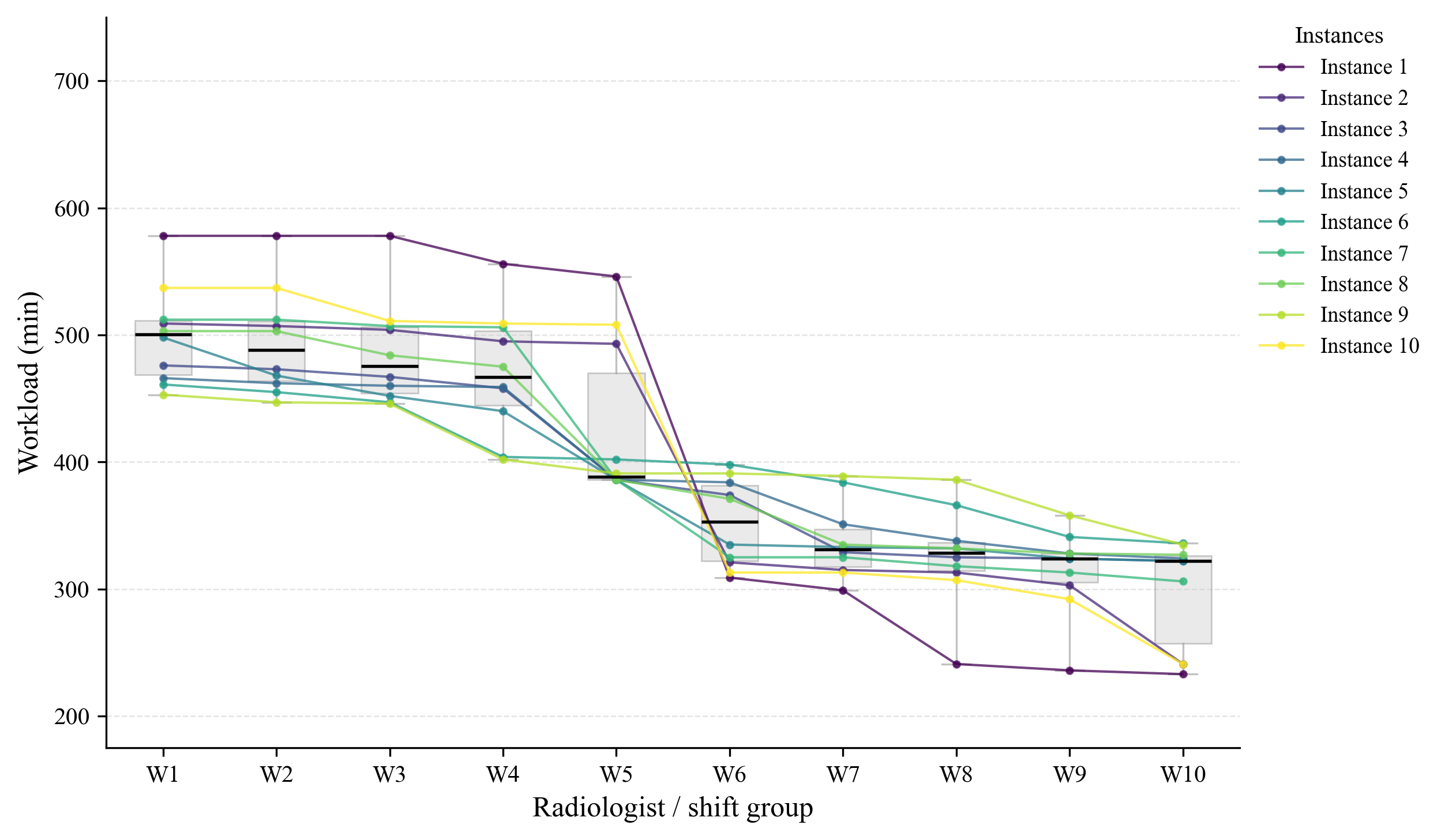}
    \caption{Workloads (in minutes) of the radiologists with $\varepsilon$ fixed as (iv).}
    \label{fig:FAIR:inst4}
\end{subfigure}

\vspace{0.5cm}

\begin{subfigure}{0.75\textwidth}
    \centering
    \includegraphics[width=\linewidth]{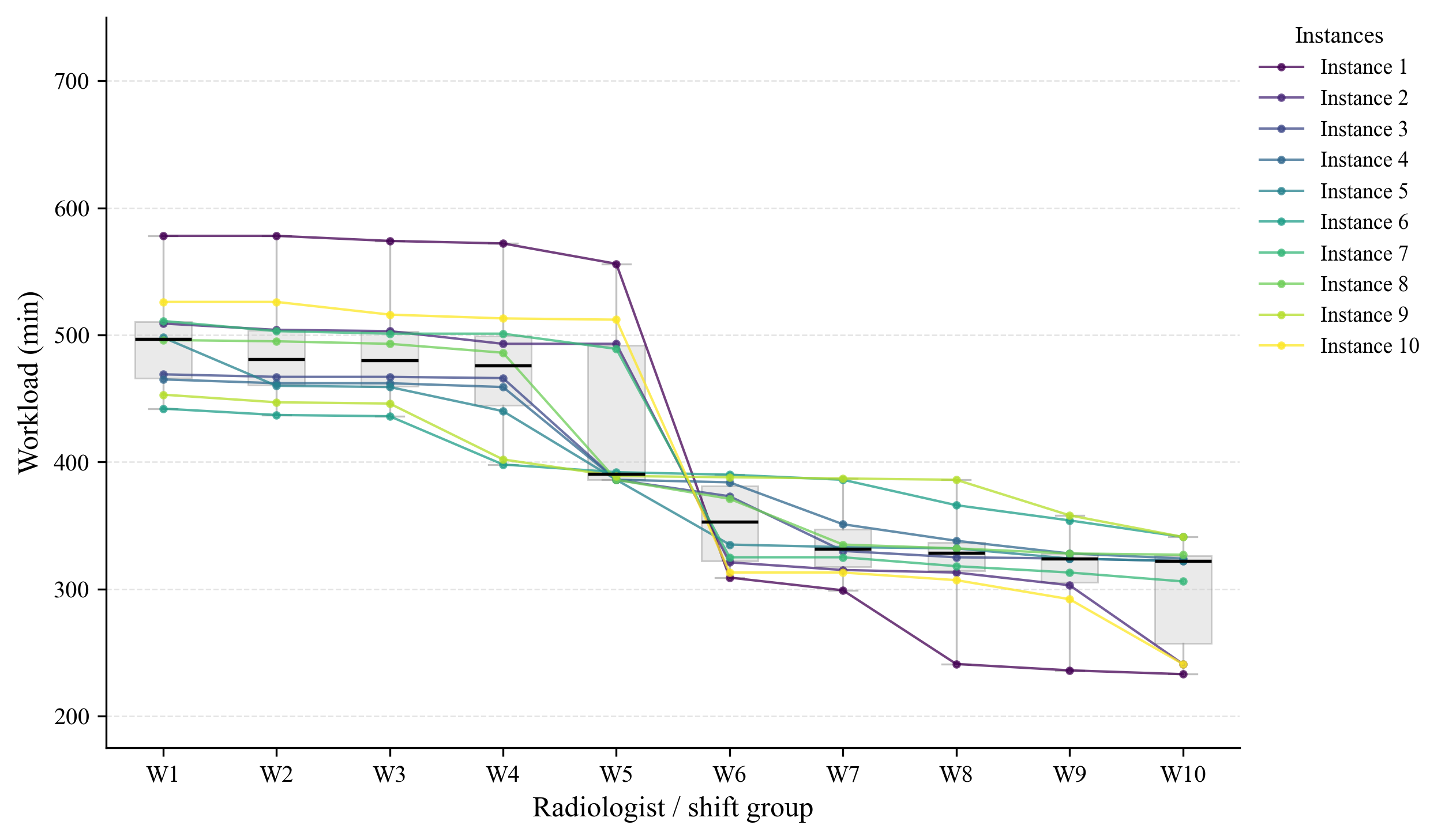}
    \caption{Workloads (in minutes) of the radiologists with $\varepsilon$ fixed as (v).}
    \label{fig:FAIR:inst5}
\end{subfigure}

\caption{Workloads of the radiologists with the baseline approach and with different values of $\varepsilon$.}
\label{fig:FAIR:spaghetti}
\end{figure}

With configuration (i), the proposed model keeps the maximum workload below 700 minutes, while the minimum workload remains above 200 minutes. Increasing the scheduling flexibility to configuration (ii) further improves workload balance: the maximum workload decreases below 600 minutes, whereas the minimum workload increases to approximately 200--300 minutes.

Configuration (iii) provides the largest improvement in workload balance. In this case, the workloads clearly separate into two clusters corresponding to the two radiology specialties. The higher-workload cluster consists of the body radiologists, who report the larger share of outpatient examinations, while the lower-workload cluster corresponds to the neuro radiologists. Within each specialty, workloads are distributed more evenly than in the previous configurations.

Allowing larger deviations from the preferred examination day, as in configurations (iv) and (v), results in only marginal additional improvements. Although slight differences can be observed across individual instances, the workload distributions remain largely comparable to those obtained with configuration (iii).

\begin{figure}
    \centering
    \includegraphics[width=0.9\linewidth]{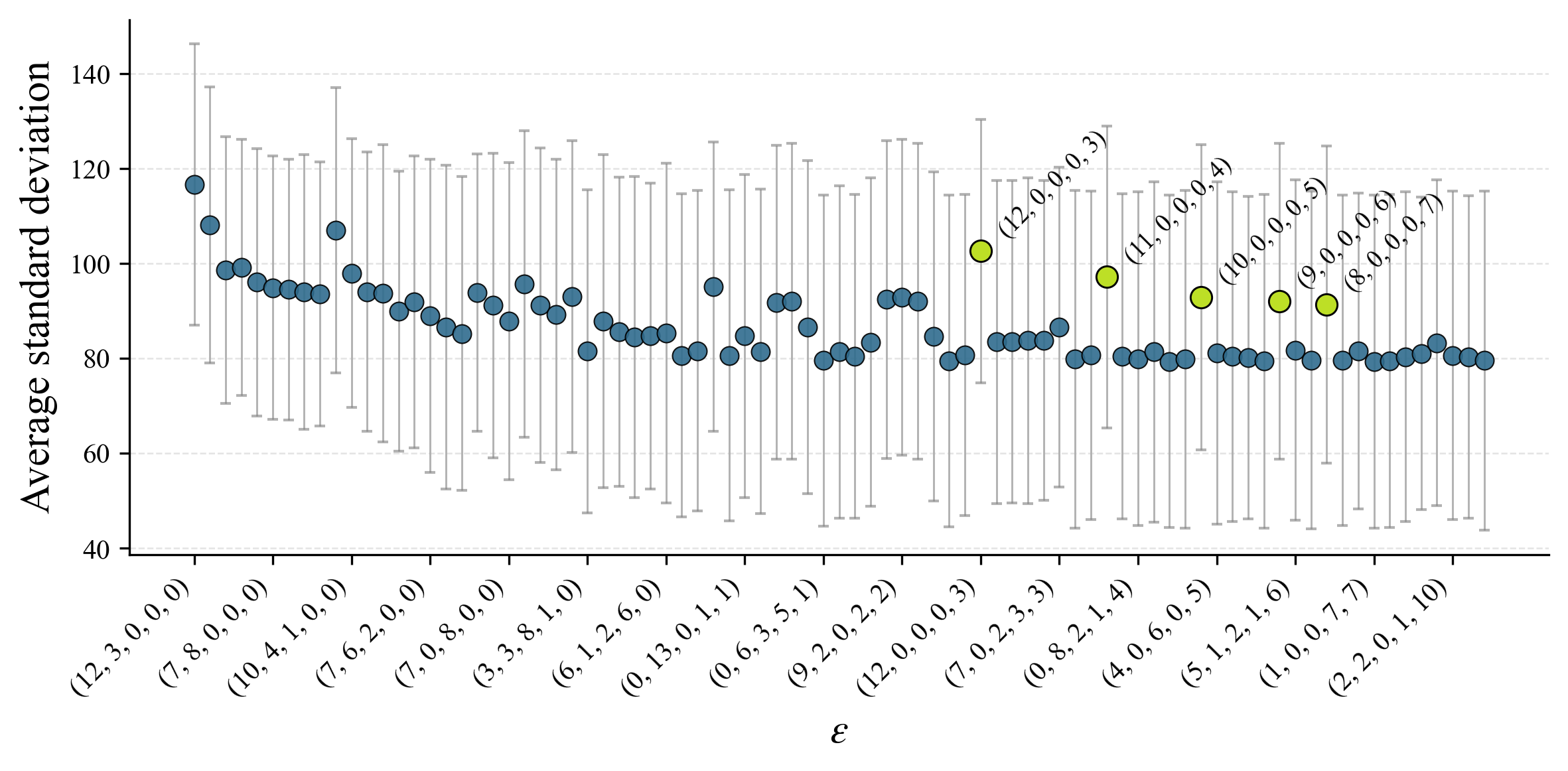}
    \caption{Average standard deviation of the radiologist workloads across the instances of the validation set as a function of  $\varepsilon$.}
    \label{fig:FAIR:pareto_std}
\end{figure}

Figure~\ref{fig:FAIR:pareto_std} reports the average standard deviation of the radiologist workloads across one instances of the test set as a function of the selected value of $\varepsilon$. Overall, larger values of $\varepsilon$ (in the lexicographic sense) are associated with lower workload standard deviations, indicating a more balanced workload distribution. A few outliers can nevertheless be observed, corresponding to $\varepsilon$ configurations that prioritize scheduling a large number of patients on their preferred examination day. In these cases, the reduced scheduling flexibility limits the model's ability to balance workloads effectively, resulting in higher workload variability.

Overall, these results indicate that the proposed optimization model effectively balances radiologists' workloads while satisfying patient--radiologist compatibility constraints. Moreover, permitting appointments to be postponed by up to two days appears sufficient to achieve most of the attainable improvement in workload fairness, whereas additional scheduling flexibility provides only limited further benefits.






\subsubsection{Scalability analysis}
We evaluate the scalability of the proposed model under increasing patient demand and resource availability. Starting from the case-study setting, which consists of 2 CT scanners, 2 radiologists (one per specialty), and 15 outpatients per week, we consider a larger scenario with 4 CT scanners, 4 radiologists (two per specialty), and a proportionally increased outpatient demand. Patients' preferred examination days are generated according to a discrete uniform distribution over the planning horizon, while the specialty distribution observed in the case study is preserved. 

\begin{figure}[tbp]
    \centering
    \includegraphics[width=0.9\textwidth]{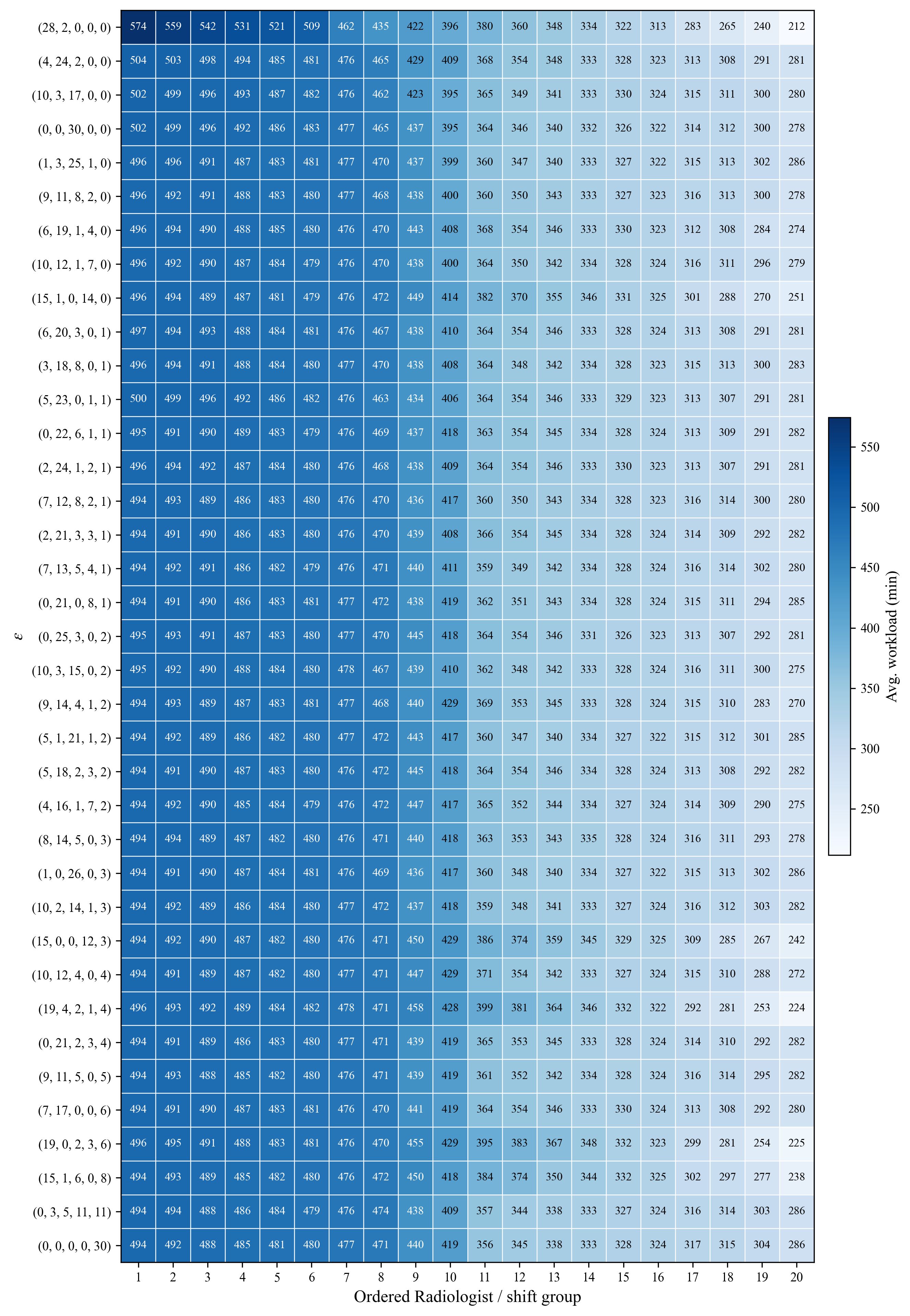}
    \caption{Heatmap of the workload (in minutes) of radiologist-shift couples, depending on $\varepsilon$, in the large instance.}
    \label{fig:FAIR:heat_double}
\end{figure}

Figure~\ref{fig:FAIR:heat_double} presents the heatmaps of the radiologist workloads for representative values of $\varepsilon$ in the larger scenario. Specifically, we report the smallest feasible value of $\varepsilon$, the largest value of $\varepsilon$, and one representative value every 50 non-Pareto-dominated $\varepsilon$ configurations.

In this larger instance, the computational time increases dramatically, with the average runtime rising from 3 minutes to 645 minutes. Indeed, the cardinality of $\Omega$ is 46,375, and thanks to the dominance reduction, an average of 29,251 MILPs can be skipped, corresponding to a reduction of approximately $63\%$ in the computational effort.

The smallest value of $\varepsilon$ for which all instances are feasible is $(28,2,0,0,0)$, indicating that 28 patients are scheduled on their preferred day, while 2 patients are shifted by one day. This represents the most unbalanced solution. 

\begin{figure}[tbp]
    \centering
    \includegraphics[width=0.9\textwidth]{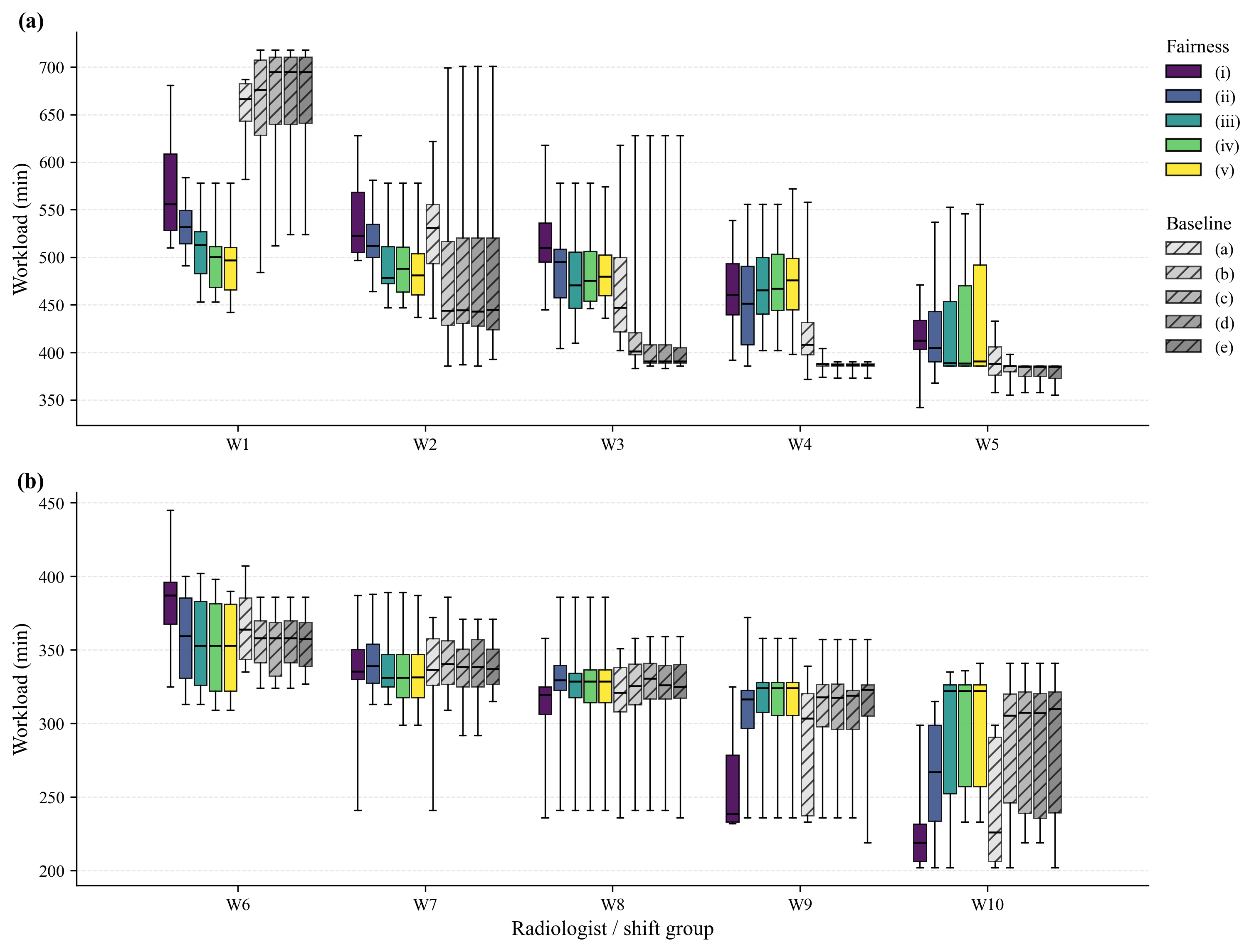}
    \caption{Workloads (in minutes) of the radiologists in the small instance set.}
    \label{fig:FAIR:boxplot_workload_small}
\end{figure}

\begin{figure}[htb]
    \centering
    \includegraphics[width=0.9\textwidth]{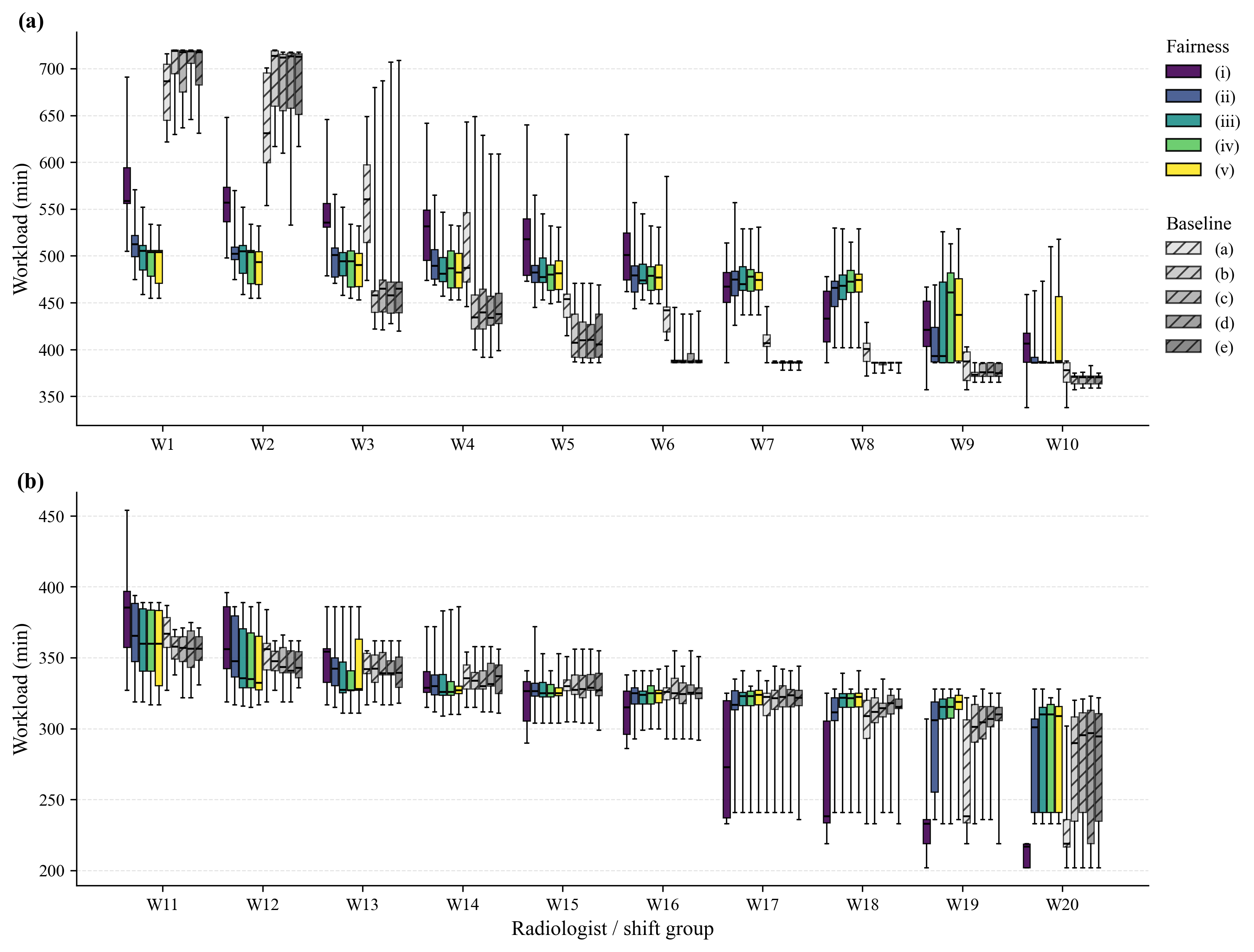}
     \caption{Workloads (in minutes) of the radiologists in the large instance.}
    \label{fig:FAIR:boxplot_workload_big}
\end{figure}

Figures~\ref{fig:FAIR:boxplot_workload_small} and~\ref{fig:FAIR:boxplot_workload_big} present the boxplots of the ordered radiologist--day workloads across the instances of the test set for the small- and large-scale scenarios, respectively. For the fairness-aware optimization model, we consider the $\varepsilon$ configurations (i), (ii), (iii), (iv), and (v). For the baseline model, we consider the following configurations: (a) the minimum total number of waiting days for which the problem is feasible (equal to 3 for the small instance and 2 for the large one); (b) the total number of waiting days equals the number of patients; (c) the total number of waiting days equals twice the number of patients; (d) the total number of waiting days equals three times the number of patients; and (e) the total number of waiting days equals four times the number of patients.

The results show that the proposed fairness-aware model consistently achieves a more balanced workload distribution than the baseline approach, particularly for body radiologists, who are responsible for the largest share of outpatient examinations. As the value of $\varepsilon$ increases, workload balance further improves because the model has greater flexibility in assigning patients to examination days. These findings indicate that the proposed approach maintains its effectiveness as the problem size increases, demonstrating its scalability to larger radiology departments within a reasonable computational time.

\subsubsection{Specialty case-mix}
We first evaluate the proposed model using the case-mix observed in the real case study, consisting of 35.9\% neuro outpatients and 64.1\% body outpatients. We then consider two alternative scenarios: a balanced case-mix (50\% neuro and 50\% body outpatients) and an unbalanced case-mix (25\% neuro and 75\% body outpatients). These scenarios are designed to assess the impact of patient--radiologist compatibility constraints on scheduling performance. In all experiments, the problem size is fixed to that of the case study, while patients' preferred examination days are generated according to a discrete uniform distribution over the five-day planning horizon.

\begin{figure}
    \centering
    \includegraphics[width=0.8\linewidth]{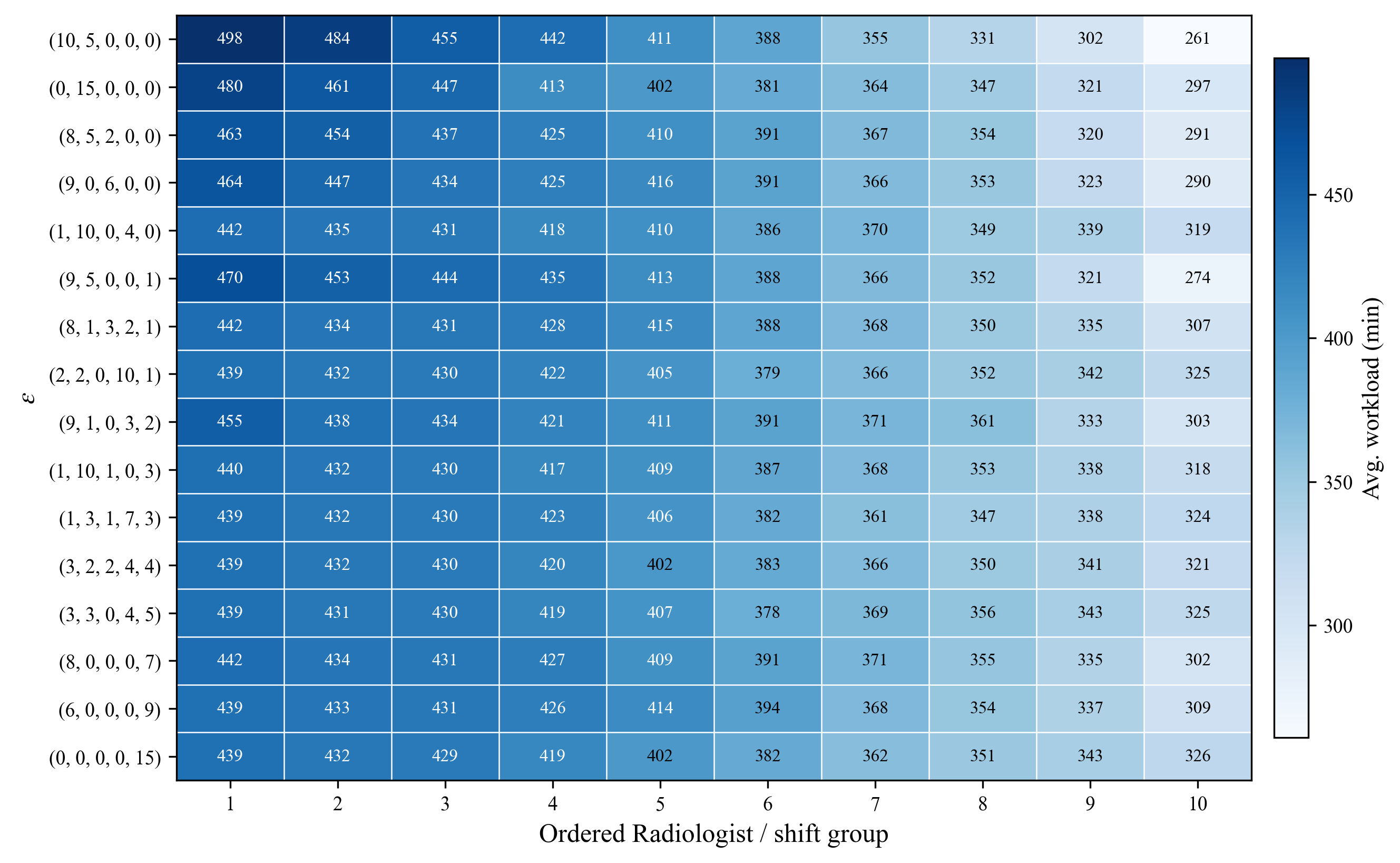}
    \caption{Heatmap of the workload (in minutes) of radiologist-shift couples, for some representative $\varepsilon$ for the balanced casemix.}
    \label{fig:FAIR:heat_casemix_50}
\end{figure}

\begin{figure}
    \centering
    \includegraphics[width=0.8\linewidth]{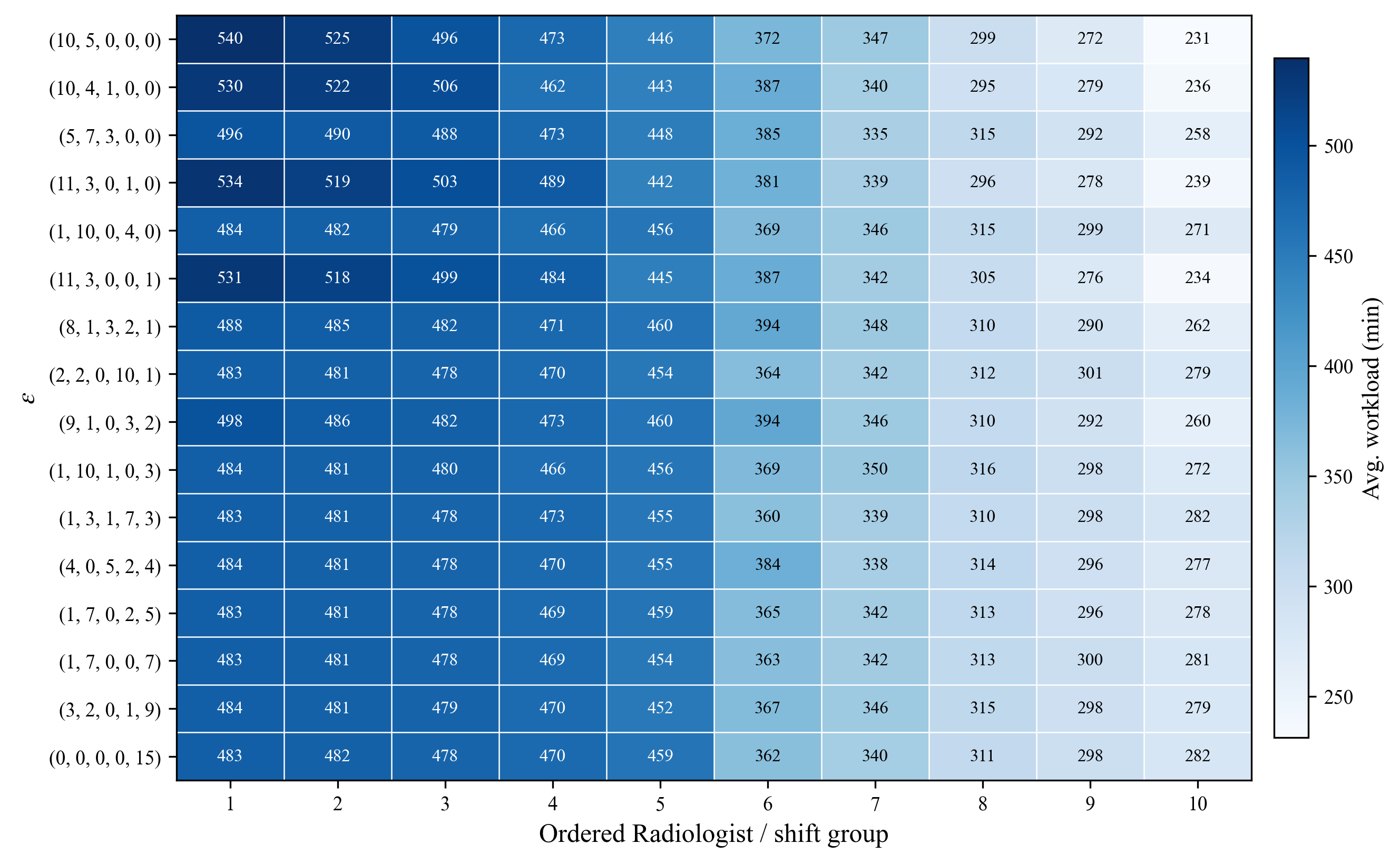}
    \caption{Heatmap of the workload (in minutes) of radiologist-shift couples, for some representative $\varepsilon$ for the unbalanced casemix.}
    \label{fig:FAIR:heat_casemix_75}
\end{figure}

Figures~\ref{fig:FAIR:heat_casemix_50} and~\ref{fig:FAIR:heat_casemix_75} present the corresponding heatmaps for a representative subset of $\varepsilon$ values. As expected, the balanced case-mix produces a more even workload distribution, since the demand is more uniformly distributed across the two radiology specialties. Conversely, increasing the proportion of body examinations naturally leads to greater workload disparities due to the patient--radiologist compatibility constraints. Nevertheless, the proposed fairness-aware approach is able to effectively mitigate these imbalances, achieving substantially more balanced workload distributions even in the most unbalanced scenario, although the attainable level of fairness is inherently limited by the specialty compatibility requirements.

\begin{figure}[htb]
    \centering
    \includegraphics[width=0.9\textwidth]{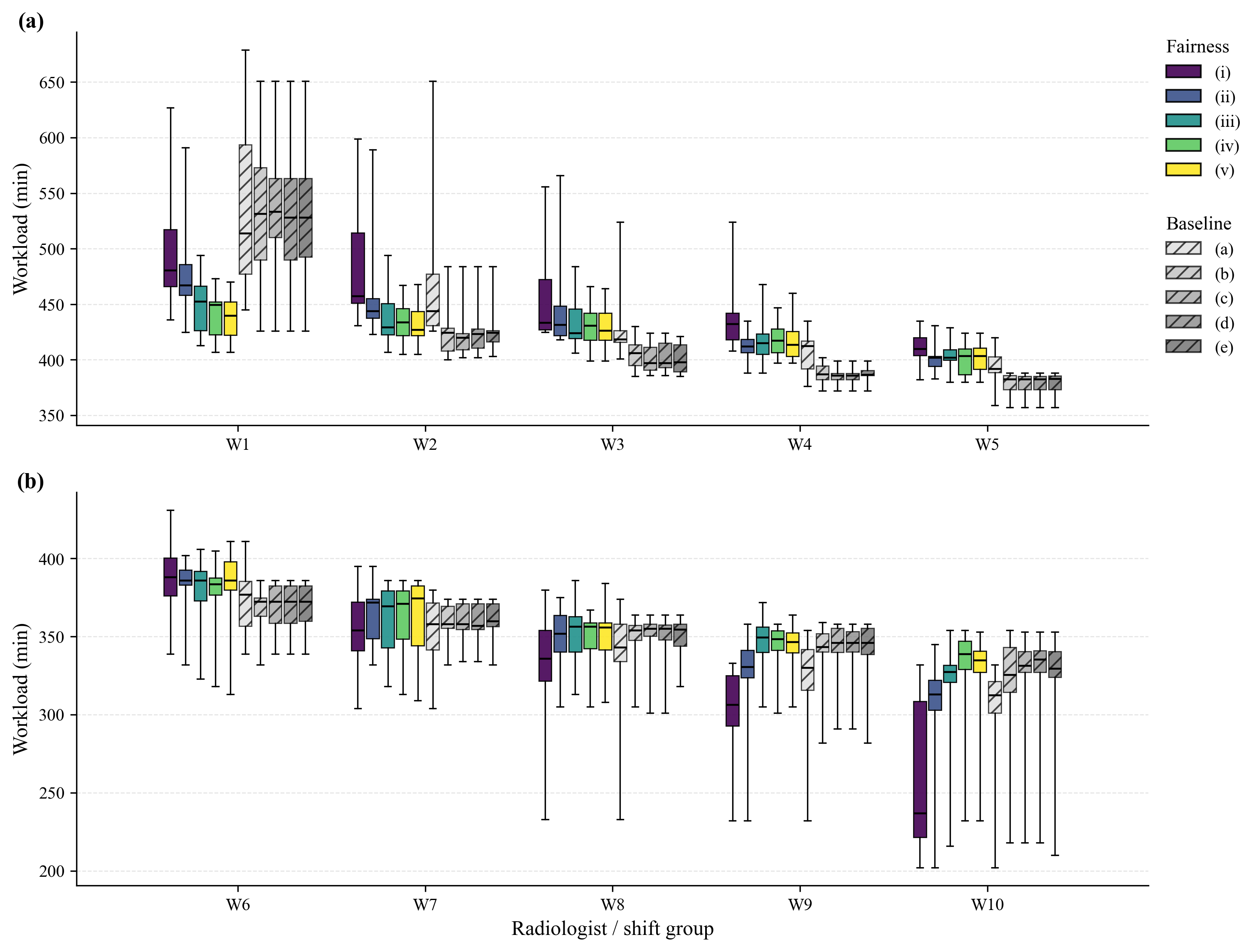}
     \caption{Workloads (in minutes) of the radiologists in the balanced case-mix.}
    \label{fig:FAIR:boxplot_casemix_50}
\end{figure}

\begin{figure}[htb]
    \centering
    \includegraphics[width=0.9\textwidth]{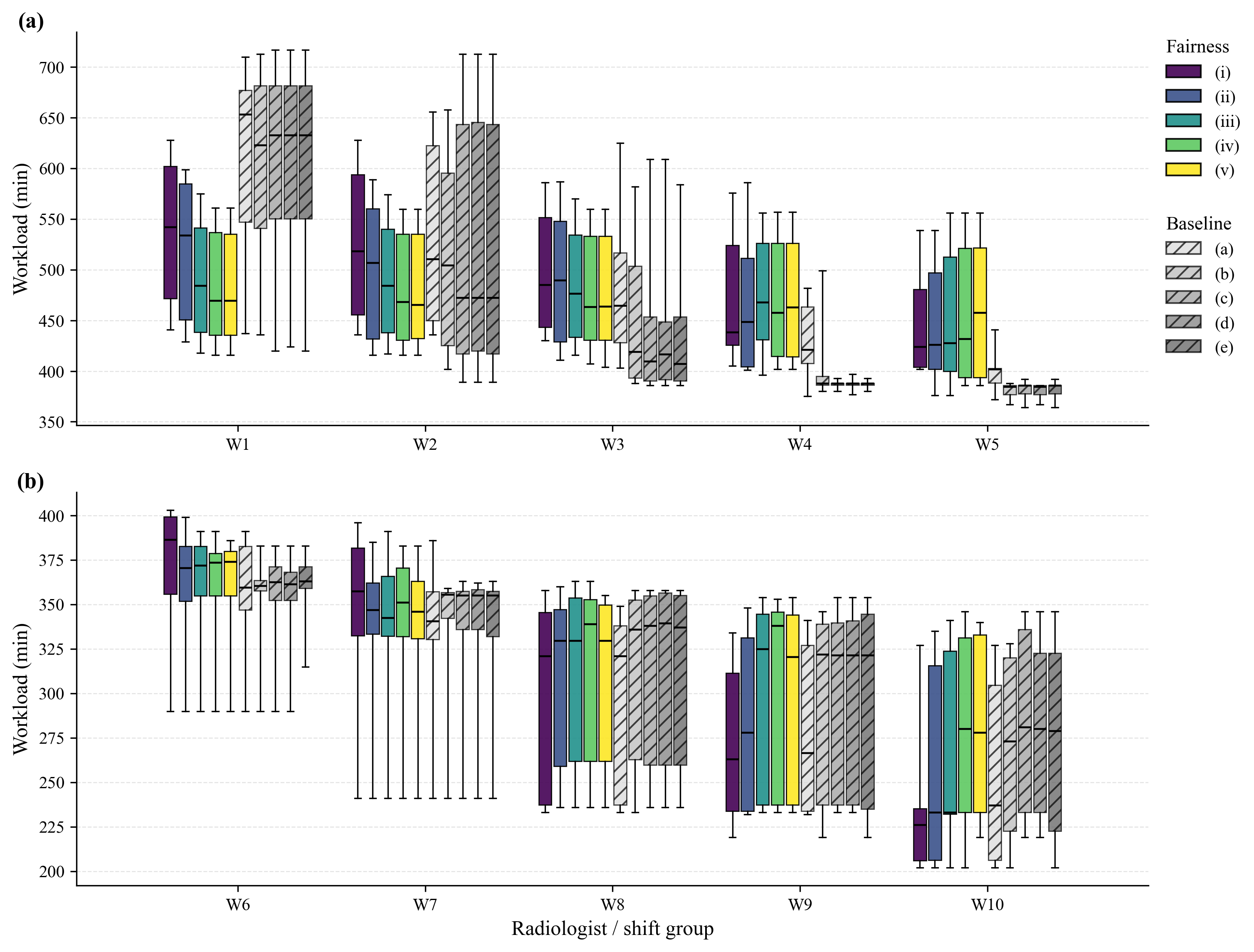}
     \caption{Workloads (in minutes) of the radiologists in the unbalanced case-mix.}
    \label{fig:FAIR:boxplot_casemix_75}
\end{figure}

 Figures~\ref{fig:FAIR:boxplot_casemix_50} and~\ref{fig:FAIR:boxplot_casemix_75} present the workload distributions obtained under the two alternative case-mix scenarios. As in the previous analysis, for the fairness-aware optimization model we consider the following $\varepsilon$ configurations: (i) $\varepsilon=(10,5,0,0,0)$, (ii), (iii), (iv), and (v). For the baseline model, we consider the corresponding configurations: (a) $\tilde{\varepsilon}=5$, (b) a total number of waiting days equal to the number of patients, (c) twice the number of patients, (d) three times the number of patients, and (e) four times the number of patients.

When the case-mix is balanced (Figure~\ref{fig:FAIR:boxplot_casemix_50}), the distinction between body and neuro radiologists effectively disappears, and no clear workload separation between the two specialties is observed. In this setting, workloads are distributed more evenly across all radiologists, typically ranging between 300 and 500 minutes, with the exception of configuration (i), which yields a less balanced workload distribution. Consequently, workload differences are no longer primarily driven by radiologist specialization.

This scenario can also be interpreted as the limiting case in which patient--radiologist compatibility constraints are absent, allowing any radiologist to report any examination.

Even under these favorable conditions, the baseline approach fails to achieve an equitable workload distribution, highlighting the effectiveness of the proposed fairness-aware optimization model.

In the unbalanced case-mix scenario, radiologist workloads remain clearly clustered according to specialty. The proposed approach effectively balances the workloads of body radiologists, whose higher workload reflects the larger proportion of body examinations. In contrast, the workload distribution among neuro radiologists remains less balanced.

A plausible explanation lies in the lexicographic optimization process. Since body radiologists experience the highest workloads, they are prioritized during the workload balancing procedure. Neuro patients must then be scheduled within the residual CT scanner capacity, which is shared across both specialties. As a result, once the workloads of body radiologists have been balanced, neuro examinations are assigned to the remaining available time slots. This reduces the scheduling flexibility available for balancing the workloads of neuro radiologists and, consequently, limits the level of fairness that can be achieved within that specialty.

\subsubsection{Preferred shift distribution}
To evaluate the impact of patient scheduling preferences, we examine three distinct probability distributions for the preferred examination day over a five-day planning horizon: a \textit{uniform} preference across all days, represented by $p_d = 0.2$ for $d \in \{1, \dots, 5\}$; A \textit{mid-week centered} preference, with $p_3 = 0.6$ and $p_d = 0.1$ for $d \neq 3$; An \textit{early-week concentrated} preference, with $p_1 = 0.6$ and $p_d = 0.1$ for $d \neq 1$.

In all computational experiments, the problem dimensions and specialty breakdown are fixed to match the real-world case study.

\begin{figure}[tbp]
    \centering
    \includegraphics[width=0.9\textwidth]{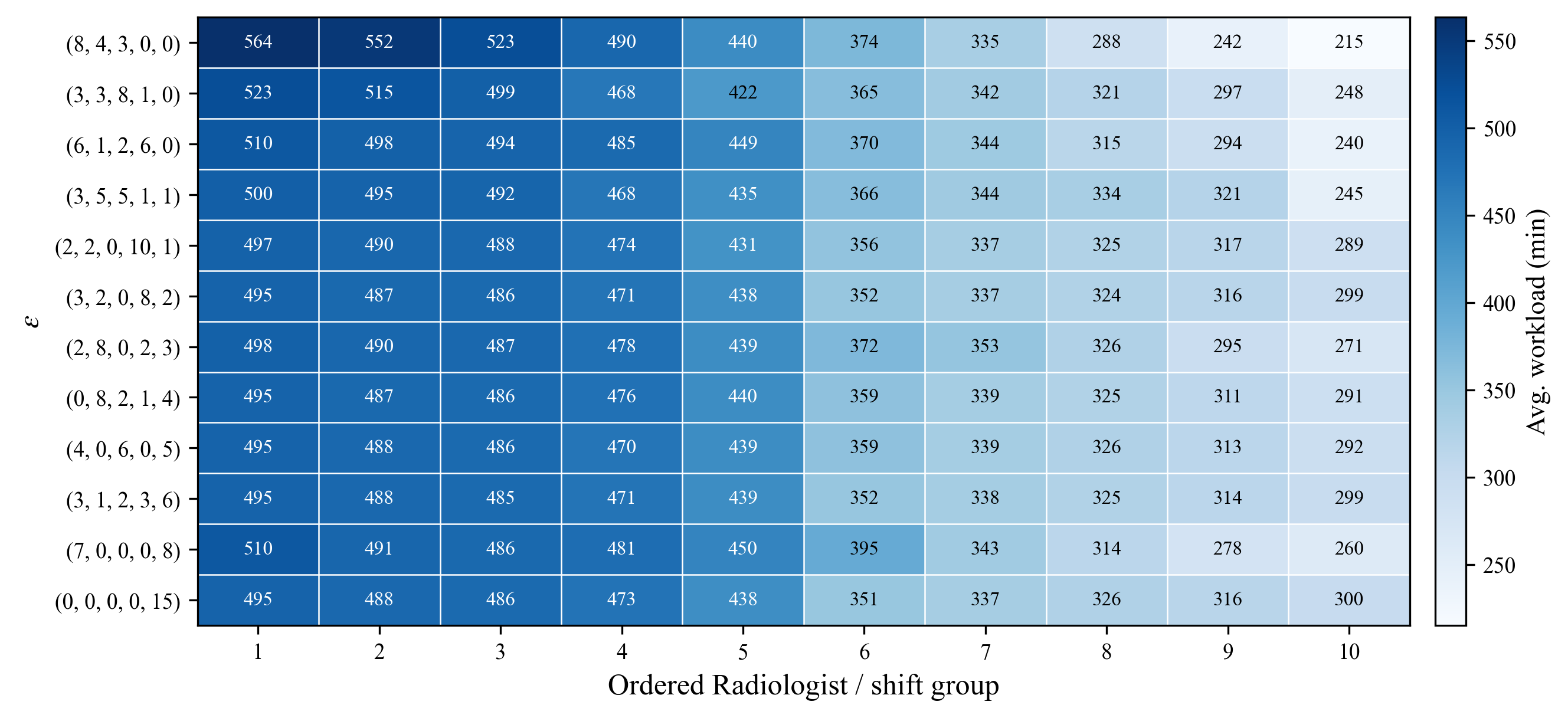}
     \caption{Heatmap of the workload (in minutes) of radiologist-shift couples, with a preferred shift distribution centered in the first day.}
    \label{fig:FAIR:heat_first}
\end{figure}

\begin{figure}[htb]
    \centering
    \includegraphics[width=\textwidth]{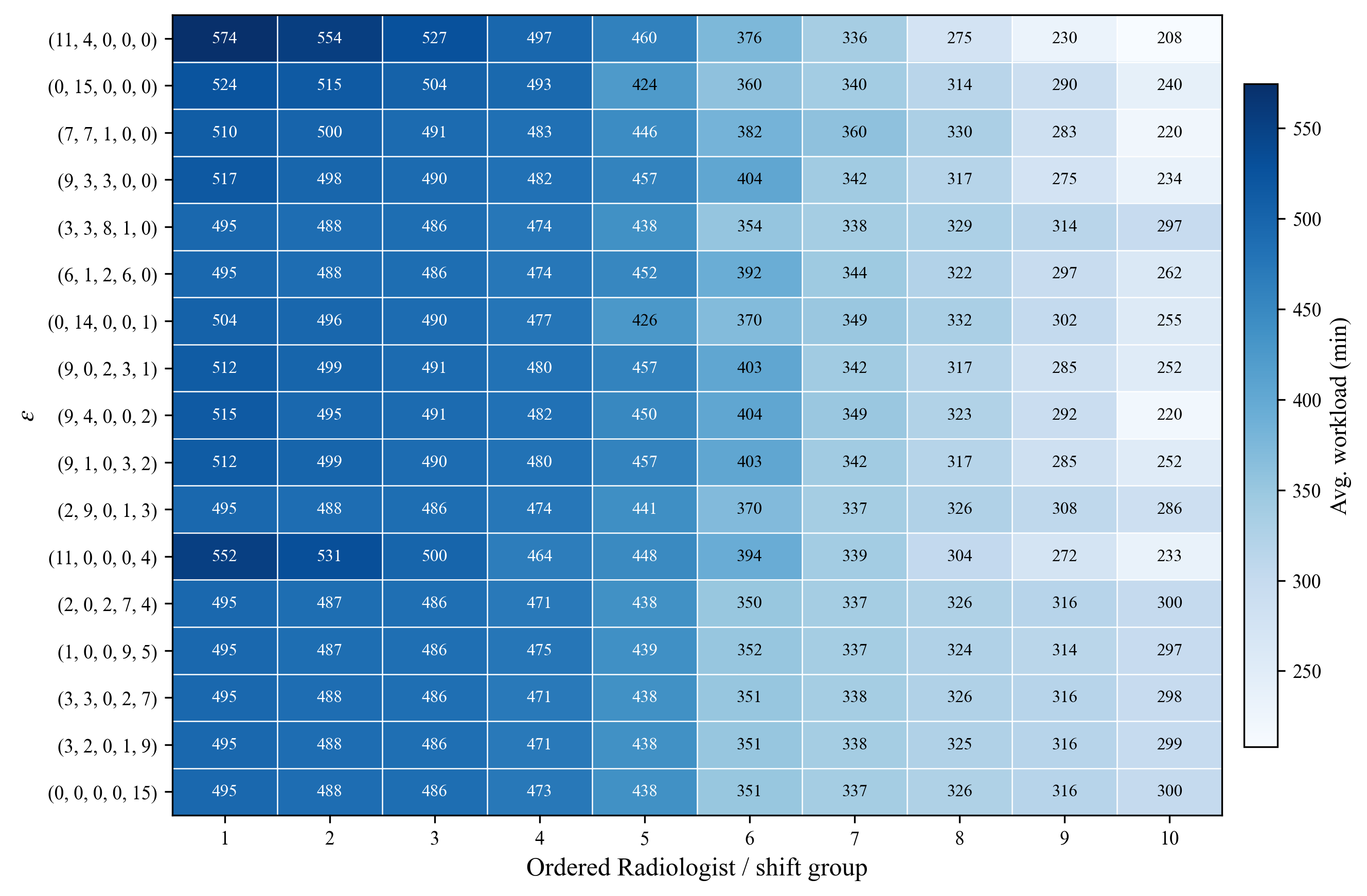}
     \caption{Heatmap of the workload (in minutes) of radiologist-shift couples, with a preferred shift distribution centered in the middle day.}
    \label{fig:FAIR:heat_middle}
\end{figure}

\begin{figure}[htb]
    \centering
    \includegraphics[width=1\textwidth]{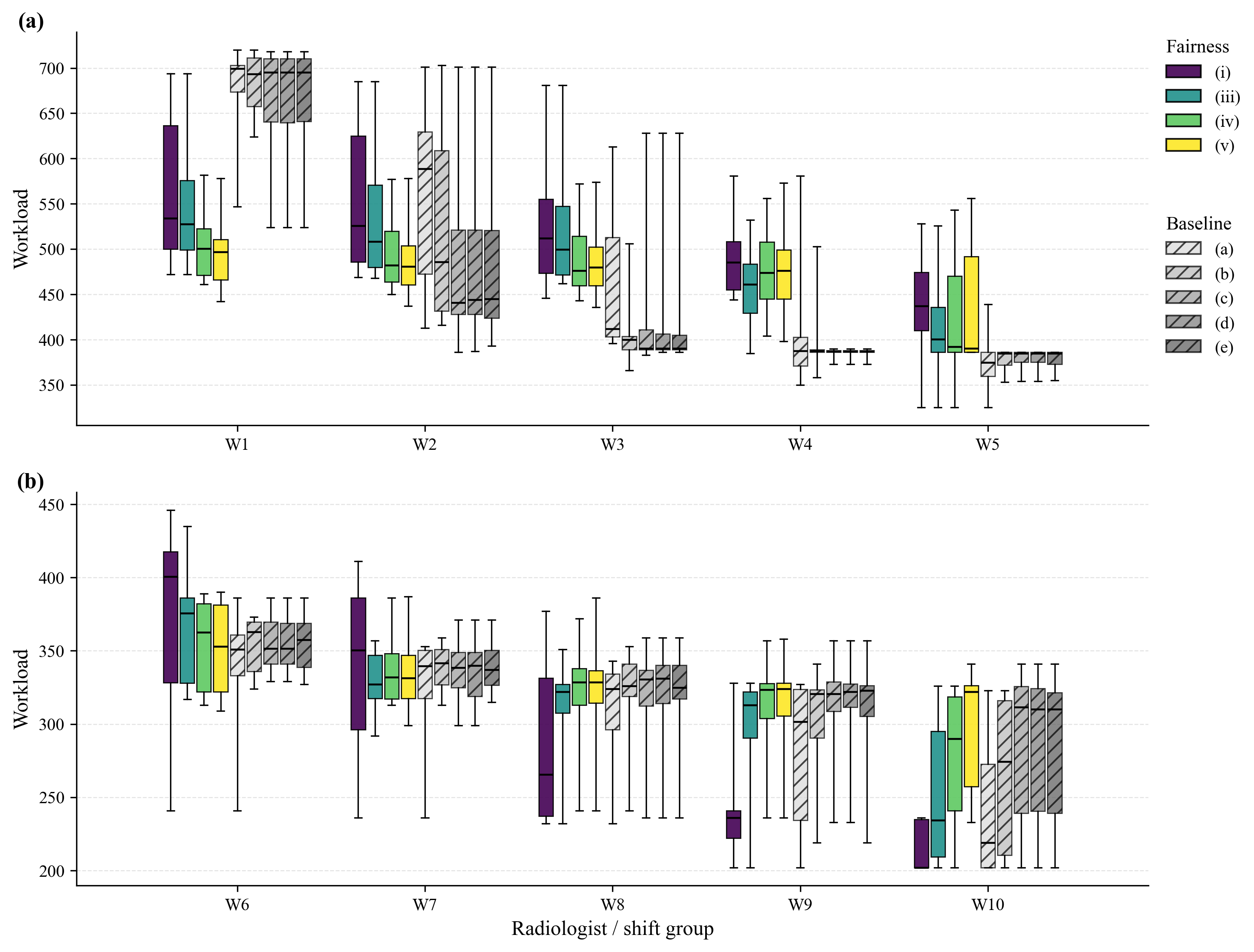}
     \caption{Workloads (in minutes) of the radiologists in the medium instance with a preferred shift distribution centered in the first day.}
    \label{fig:FAIR:boxplot_first}
\end{figure}

\begin{figure}[htb]
    \centering
    \includegraphics[width=0.9\textwidth]{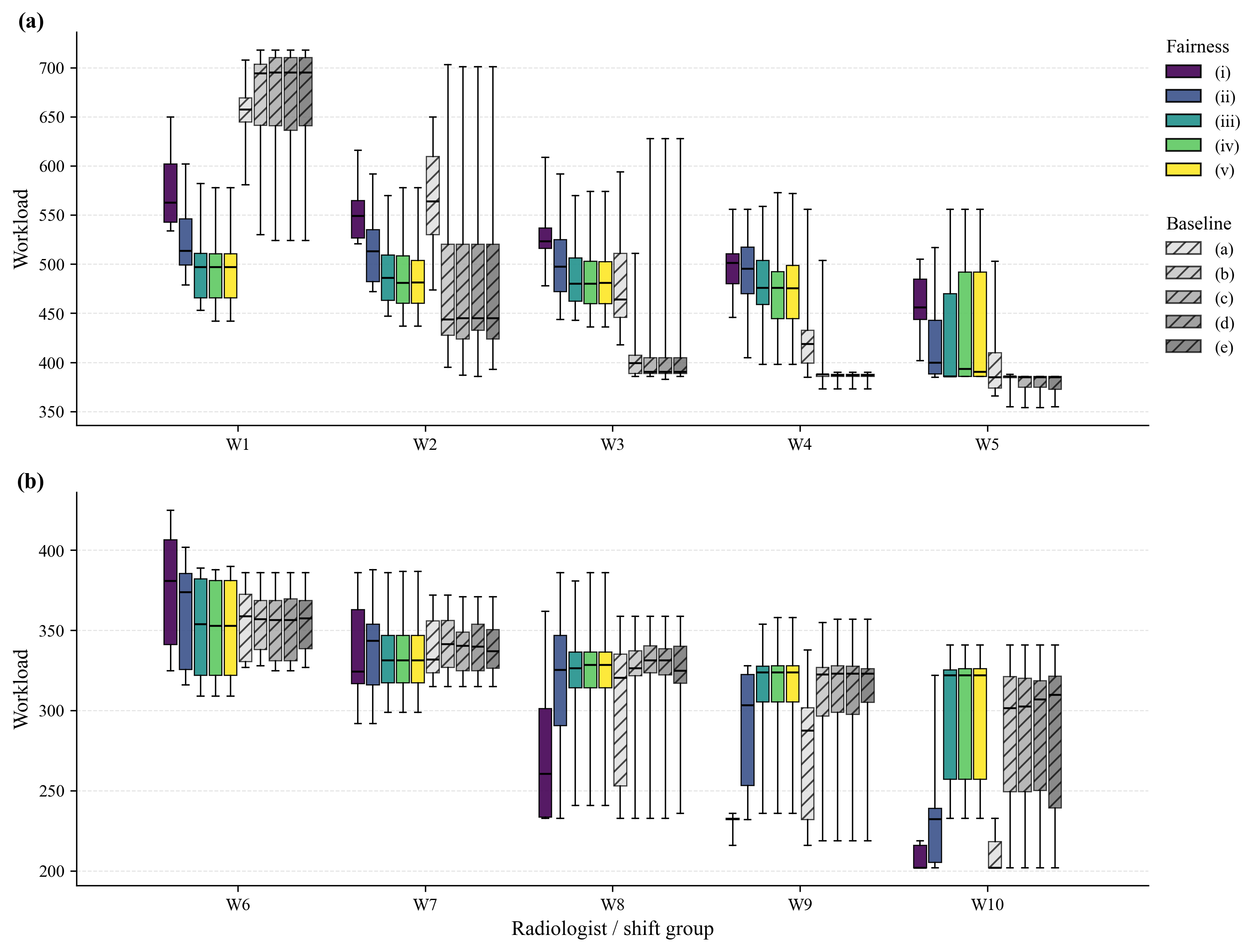}
     \caption{Workloads (in minutes) of the radiologists in the medium instance with a preferred shift distribution centered in the middle day.}
    \label{fig:FAIR:boxplot_middle}
\end{figure}

The results for the uniform preference distribution are reported in Figures~\ref{fig:FAIR:heat_case_study} and~\ref{fig:FAIR:boxplot_workload_small}. The corresponding heatmaps for the two alternative preference distributions are shown in Figures~\ref{fig:FAIR:heat_first} and~\ref{fig:FAIR:heat_middle}, while the associated workload boxplots are presented in Figures~\ref{fig:FAIR:boxplot_first} and~\ref{fig:FAIR:boxplot_middle}.

In this setting, although the maximum workload remains relatively stable, larger values of $\varepsilon$ are required to obtain feasible configurations. This indicates that, to achieve a comparable level of workload balancing, greater flexibility in patient scheduling is necessary. The two configurations produce broadly similar workload distributions; however, as expected, when patients' preferences are concentrated on the first day, a larger number of patients must be shifted from their preferred examination day to obtain a feasible solution.

When the preferred examination day distribution is concentrated on the first day, no feasible solution is obtained for configuration (ii), while configuration (i) requires $\varepsilon=(8,4,3,0,0)$ to achieve feasibility. As $\varepsilon$ increases, a substantial improvement in workload balance is observed. This behavior differs from the uniform preference distribution, where configurations (iv) and (v) provided results that were largely comparable to those obtained with configuration (iii).

When the preferred appointment-day distribution is concentrated around the central day, configuration (i) corresponds to $\varepsilon=(11,4,0,0,0)$, and feasible solutions can already be obtained for configuration (ii). However, the quality of these solutions is significantly lower than that achieved with larger values of $\varepsilon$. The remaining $\varepsilon$ configurations exhibit comparable performance: configuration (iii) provides slightly inferior results, whereas configurations (iv) and (v) show almost identical behavior.

\subsection{Managerial Insights}
The integration of predictive models into the multi-objective scheduling framework yields several important operational and managerial insights for hospital administrators.

From the predictive perspective, the first key finding is that evaluating ML models solely through traditional performance metrics (e.g., $\text{MAE}$ or $R^2$) can be misleading in prescriptive applications. A model with a low mean absolute error may still exhibit systematic directional bias by consistently overestimating or underestimating process durations. In a predict-then-optimize framework, systematic overestimation creates artificial bottlenecks and increases scanner idle time, whereas systematic underestimation leads to schedule disruptions, patient delays, and staff overtime. XGBoost achieves the lowest decision regret precisely because it minimizes schedule-level operational disruptions rather than unweighted prediction errors, making it the most reliable model for downstream decision-making.

From the prescriptive perspective, the optimization model effectively highlights the trade-off between patient preference satisfaction and radiologist workload equity. Computational results show that allowing a modest appointment flexibility of only 1-2 days from the patient's preferred examination date is sufficient to achieve a balanced workload among radiologists. This finding suggests that a high degree of workload equity can be attained without requiring substantial deviations from patient preferences.

In the analyzed case study, the primary operational challenge arises from workload imbalances across radiological specialties. Since patients requiring body CT examinations constitute the majority of the demand, the optimization model prioritizes balancing the workload of body radiologists. As a result, the remaining scheduling flexibility for neuroradiologists is reduced, naturally limiting the degree of workload equity that can be achieved within this specialty. Therefore, hospital managers should recognize that the demand distribution across clinical specialties determines the priority and effectiveness of workload balancing.

When patient preferences or the clinical case-mix are strongly concentrated on specific days or specialties, the optimization model explicitly reveals the resulting structural workload asymmetries. If these demand patterns persist over time, the proposed framework can serve as a valuable diagnostic tool for tactical capacity planning, enabling decision-makers to assess whether workforce allocation, radiologist shift planning, or scanner time allocation should be strategically adjusted to better align capacity with long-term demand trends.

Although analyzing all possible combinations of allowable waiting times requires solving a very large number of MILP instances, the proposed dominance reduction property significantly improves computational efficiency. Realistic problem instances are solved rapidly, and even under stress-testing conditions, the solution time remains well within acceptable operational limits (645 minutes). Since CT master schedules are typically generated on a weekly basis, this computational performance supports practical real-world implementation and enables efficient what-if and scenario analyses.
\section{Conclusions}
\label{sec:FAIR:conclusion}

The proposed approach provides decision-makers with a comprehensive overview of the operational consequences of different scheduling policies. By integrating ML-based predictions with the optimization framework, the model enables the assessment of the trade-off between patient flexibility and radiologist workload balancing.

The framework supports both short-term and long-term managerial decisions. At the operational level, it helps managers identify the most appropriate weekly scheduling strategy by evaluating how different levels of patient scheduling flexibility affect workload distribution and waiting times. At the strategic level, it provides insights into structural workload imbalances arising from factors such as specialty case-mix or concentrated patient preferences, supporting decisions related to workforce allocation and capacity planning.

The analysis of different scenarios shows that a more balanced case-mix naturally leads to more homogeneous workload distributions, while highly concentrated patient preferences may require greater scheduling flexibility to achieve comparable workload balance. Therefore, the proposed framework can also serve as a diagnostic tool to identify potential sources of inefficiency and guide future organizational improvements.

Unlike approaches that provide a single scheduling solution, the proposed framework allows decision-makers to explore the trade-offs between competing objectives, increasing transparency and facilitating the adoption of optimization-based recommendations in clinical environments.

As future research directions, we aim to extend the model by introducing patient-specific delay penalties that account for the clinical impact and acceptability of postponing examinations by different amounts of time. This would allow the framework to better capture the heterogeneous consequences of waiting times across patients.

Furthermore, a user-friendly decision-support tool could be developed to be routinely used by radiology managers to automatically analyze the current operational status of the system. Beyond providing the optimized schedule, such a tool would offer insights into workload distribution, demand patterns, and potential sources of imbalance, supporting proactive decision-making and continuous improvement of radiology operations.

\section*{Acknowledgments}
The work of the author Sara Cambiaghi is supported by the Ph.D.\ fellowship funded under the Italian National Recovery and Resilience Plan (PNRR), as part of the program for research doctorates and innovative doctorates for public administration and cultural heritage (D.M.\ 118, 02-03-2023).

\singlespacing
\bibliographystyle{elsarticle-num}
\bibliography{biblio}

\end{document}